\documentclass[11pt]{article}

\usepackage{amsmath}
\usepackage{amsfonts}
\usepackage{amsthm}
\usepackage{amsbsy}

\usepackage[english]{babel} 

\usepackage{amssymb}
\usepackage{graphicx}     
\usepackage{cite}
\usepackage{epstopdf}
\usepackage{epsfig}
\usepackage{subfigure}
\usepackage{smartdiagram}

\def\Co{{\rm Cost}}

\def\E{\mathbb{E}}
\def\EE{\mathbb{E}}
\def\var{{\rm Var}}
\def\cov{{\rm Cov}}

\def\LL{{\mathcal{B}(\mathcal{R})}}

\def\llambda{\Lambda}

\def\w{v}
\def\theta{\red{z}}
\def\G{\Psi}

\def\LB{{\mathcal{B}(\mathcal{R};L_2(\Omega))}}

\newcommand{\bea}{\begin{eqnarray}}
\newcommand{\eea}{\end{eqnarray}}

\newcommand{\beas}{\begin{eqnarray*}}
\newcommand{\eeas}{\end{eqnarray*}}

\def\e{\varepsilon}

\def\cir{\makebox[0.5cm]{$\circ$}}
\def\ind{\hspace{0.25in}}

\def\nt{n_{\mbox{\sc \small tot}}}

\def\fP{\em \sf}
\def\IR{\mathop{\mbox{\rm Iround}}\nolimits}

\newcommand{\beq}{\begin{eqnarray*}}
\newcommand{\eeq}{\end{eqnarray*}}
\newcommand{\bd}{\begin{description}}
\newcommand{\ed}{\end{description}}

\newcommand{\be}{\begin{equation}}
\newcommand{\ee}{\end{equation}}
\newcommand{\ba}{\begin{array}}
\newcommand{\ea}{\end{array}}
\renewcommand{\Sigma}{{{\sigma}}}

\newtheorem{example}{Example}
\newtheorem{proposition}{Proposition}
\newtheorem{remark}{Remark}

\newcommand{\blue}[1]{{#1}}

\newcommand{\red}[1]{{#1}}
\newcommand{\green}[1]{{#1}}

\def\beq#1{\blue{\[#1\]}}

\def\R{\mathbb R}
\def\RR{\mathbb R}

\def\ba#1#2\ea{\begin{eqnarray}#1#2\end{eqnarray}}
\def\be#1#2\ee{\begin{eqnarray}#1#2\end{eqnarray}}

\newtheorem{algorithm}{Algorithm}

\begin{document}

\title{An introduction to uncertainty quantification for kinetic equations and related problems}
\author{Lorenzo Pareschi
\thanks{Department of Mathematics and Computer Science, University of Ferrara, Via Machiavelli 30, 44121 Ferrara, Italy, {\tt lorenzo.pareschi@unife.it}}}
%
%
\maketitle

\begin{abstract}
We overview some recent results in the field of uncertainty quantification for kinetic equations and related problems with random inputs. Uncertainties may be due to various reasons, such as lack of knowledge on the microscopic interaction details or incomplete information at the boundaries or on the initial data. These uncertainties contribute to the curse of dimensionality and the development of efficient numerical methods is a challenge. After a brief introduction on the main numerical techniques for uncertainty quantification in partial differential equations, we focus our survey on some of the recent progress on multi-fidelity methods and stochastic Galerkin methods for kinetic equations.


\end{abstract}
\tableofcontents

\section{Introduction}
Many physical, biological, social, economic, financial systems involve
{uncertainties} that must be taken into account
in the mathematical models, for example partial differential equations (PDEs),
describing these systems \cite{BLMS,GWZ,JinPareschi,NPT,PIN_book, PT2, Xu}. These may be due to incomplete knowledge of the system ({epistemic uncertainties}) or they may be intrinsic to the system and cannot be reduced through improvements in measurements, etc. ({aleatoric uncertainties}). Examples include uncertainty in the \blue{initial data} and in the \blue{boundary conditions}, or in the \blue{modeling parameters}, like microscopic interactions, source terms and external forces. 
From the point of view of numerical methods, there are no relevant differences between the types and sources of uncertainty, so we will not be concerned about the nature of the uncertainties in the description of the system.

In this context, a particularly challenging case is represented by kinetic equations with random inputs. The construction of numerical methods for uncertainty quantification (UQ) in kinetic equations is a problem of considerable interest that has recently attracted the attention of many researchers (see \cite{albi2015MPE, CPZ, CZ, DJL, DPe, DPms, DPms2, DPZ, DPZ2, GJL, HLPR, HJ, JSh, JZ, LJ, LW, LX, Po, Po2, RHS} and the collection \cite{JinPareschi}). Some of the main difficulties that characterize the development of efficient methods for these equations concern the high dimensionality of the problems that, besides the variables that characterize the phase space, contain stochastic parameters, and the constraints imposed by physical properties such as the positivity of the solution and the equilibrium states. The latter problem is closely related to the construction of stochastic asymptotic preserving methods \cite{AF, JL2, JPH, JPH2, LL, ZJ}. 

In this survey, we will address some recent developments in this direction based on the use of Monte Carlo-type (MC) techniques \cite{DPms, DPms2, HPW} and on the use of Stochastic Galerkin-type (SG) approaches \cite{DPZ, DPZ2, CPZ, PZ2}. This short survey 
and the selected bibliography are obviously biased by the personal contributions and knowledge of the author, and are not intended to provide a complete overview of all the different techniques that have been developed for the quantification of uncertainty in kinetic equations. Uncertainty quantification is such a broad and active field of research that it is impossible to give credit to all relevant contributions.

The rest of the manuscript is organized as follows.
After a general introduction, in Section 2, on uncertainty quantification for PDEs and  related numerical techniques, we will focus our survey on the case of kinetic equations. In the first part, we will introduce multi-fidelity techniques to accelerate the convergence of MC methods. These techniques are particularly effective in the context of kinetic equations thanks to the presence in the literature of several surrogate models constructed with the aim of reducing the computational cost of the full model, typically represented by a Boltzmann-type collision equation. Sections 3, 4 and 5 are dedicated to these aspects. In the second part, we will address the problem of the loss of structural properties of numerical schemes in the case of intrusive SG approaches. In this context, we will first illustrate a technique based on micro-macro decomposition that allows us to efficiently and accurately approximate the equilibrium states. Subsequently, we will introduce a novel hybrid approach based on a random space SG method combined with a particle approximation of the kinetic equation in the physical space. This latter technique, allows to build efficient solvers that retain all the main physical properties, including the non-negativity of the solution.

\section{Uncertainty quantification for PDEs}
\begin{figure}[tb]
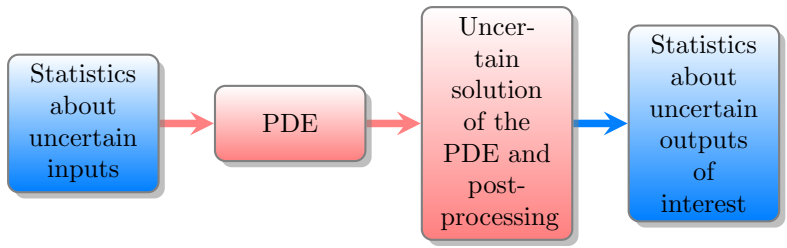

\begin{center}
\smartdiagramset{set color list={blue!50!cyan,red!50!white,red!50!white,blue!50!cyan},
   back arrow disabled=true}
\smartdiagram[flow diagram:horizontal]{{Statistics about uncertain inputs}, PDE, Uncertain solution of the PDE and post-processing,Statistics about uncertain outputs of interest}
\end{center}
\caption{The uncertainty quantification process for PDEs.}
\label{fg:fig1}
\end{figure}

The recent growth of interest in UQ for PDEs can be traced back mainly to three factors: widespread availability of data resulting from advances in technology, the increased  development of high-performance computing and the construction and analysis of new algorithms for solving differential equations with random inputs. In presence of uncertainties it becomes necessary to quantify these effects on the solution of the PDE, or on any {\em quantity of interest} (a quantity that depends on the solution of the PDE for which we want to know some statistical information), derived from the solution. The complete UQ task then consists of determining information about the uncertainty in an output of interest that depends on the solution of a PDE, given
information about the uncertainty in the inputs of the PDE (see Figure \ref{fg:fig1}).


%
 

\subsection{PDEs with random inputs}

Assume a set of \emph{random parameters} (a finite set of random numbers) $z_1,\ldots,z_{d_z}$ which may be collected in a {vector} $\red{z=(z_1,\ldots,z_{d_z})^T \in \Omega \subseteq \R^{d_z}}$.  
The solution of the PDE is not only function of the physical variables in the phase space but also of the random vector. 

For example, the scalar conservation law with random inputs
\be
\begin{aligned}
&\partial_t u(\red{z},x,t) + \nabla_x\cdot {\cal F}({\red z}, u(\red{z},x,t))= 0,
\end{aligned}
\ee
or the Fokker-Planck equation with uncertainty
\be
\begin{aligned}
\partial_t u(\red{z},x,t) + \nabla_x\cdot \left[ {\cal B}({\red z},x,t) u(\red{z},x,t)-\nabla_x \left({D}({\red z},x,t) u(\red{z},x,t)\right) \right] = 0
\end{aligned}
\ee
given the (eventually uncertain) initial data $u(\red{z},x,0)=u_0(\red{z},x)$, $x\in\R^{d_x}$, and where the terms $\cal F$, $\cal B$ and $D$ depend on the random parameters. 



A \blue{realization} of a solution of the PDE is a solution obtained for a
specific choice of the random parameters. 
One instead wants to obtain \blue{statistical information} on a quantity of interest
e.g., expected values, variances, standard deviations, covariances, higher statistical moments, etc. Therefore, multiple solutions of the PDE are necessary in order to achieve such information.

The statistical {quantities of interest} are usually determined as follows: 
\begin{itemize}
\item From the solution $u(\red{z},x,t)$ of the PDE, define an output of interest $F(u)$. 
\item Choose what statistical information about $F(u)$ is desired and define $\G(w)$ such that the {quantity of interest} is given by
\be
\E[\G(F(u))]:=\int_{\Omega} \G(F(u))p(\red{z})\,d\red{z},
\ee
where $p(\red{z})$ is the {probability density function} (PDF) of the input parameters. 
\end{itemize}
Below we give some examples of quantities of interests.
\begin{example}~
\begin{enumerate}
\item[(i)] $$F(u)=\|u\|_{L^p(\R^{d_x})}=\left(\int_{\R^{d_x}} |u(\theta,x,t)|^p\,dx\right)^{1/p},\qquad \G(w)=w$$ will give as quantity of interest $\E[\|u\|_{L^p(\R^{d_x})}]$ the {expected value} of the $L^p$-norm.  

\item[(ii)] $$F(u)=u,\qquad \G(w)=(w-\E[w])^2$$ yields $\var[u]$ the {variance} of the solution.  

\item[(iii)] $$\G(w,v)=(w-\E[w])(v-\E(v))$$ applied to $F(u_1)$ and $F(u_2)$, where $u_1$ and $u_2$ are two solutions of the PDE gives $\cov(F(u_1),F(u_2))$ the \emph{covariance}.
\end{enumerate}
\end{example}

One of the main challenges for numerical methods, is that the computational cost
associated with UQ increases with the number
of parameters used to model the uncertainty (\emph{curse of dimensionality}). This is a general problem but it is particularly relevant for \blue{kinetic equations} where the dimension of the phase-space is very high. 
We can follow two main strategies to alleviate this problem:
\begin{itemize}
\item one can try to use relatively few solutions of the PDE and replace the PDE with a surrogate, \emph{low-fidelity}, model which is much cheaper to solve. Correlation between the two models may then be used in a control variate setting.
\item for smooth solutions one can design methods which permit an accurate evaluation of $\E[\G(F(u))]$ using few quadrature points obtained from \emph{stochastic orthogonal polynomials} with respect to the PDF.  
\end{itemize}
We have tacitly assumed that we know $p(z)$, the PDF  of the
input parameters. In practice, one usually does not know much about the statistics
of the input variables and need to deal with the corresponding \emph{stochastic inverse problem} \cite{RKA}.


\subsection{Overview of techniques}

Many methods have been devised in the literature for approximating statistics of quantities of interest. We summarize shortly some of the main methods below (see \cite{BLMS, Giles, GWZ, GS, JinPareschi, LeMK, PIN_book, Xu} for recent monographs and surveys)
\begin{itemize}
\item \emph{Monte Carlo sampling}: one generates independent realizations of random inputs
based on their PDF (which may be known or not and not necessarily smooth). For each realization the problem is deterministic and can be solved by standard methods in a non intrusive way. The advantage is its simplicity but on the other hand it implies a \blue{slow convergence} and \blue{fluctuations} in the solution statistics \cite{Caflisch, HH}.  

\item \emph{Multi-fidelity, multi-level methods}: Accelerate Monte Carlo sampling methods by using \blue{multiple surrogate models} with different levels of fidelity in a control variate setting \cite{DPms, DPms2, LX, PGW, PWG}. Low-fidelity models may also be obtained from a multi-level hierarchy of numerical discretizations in the phase space \cite{FDKI, Giles, HPW, MSS}. 

\item \emph{Stochastic-Galerkin}: solutions are expressed as orthogonal polynomials of the random inputs accordingly to their PDF. Spectral convergence for smooth solutions in the random space \cite{LJ, Xu}. They require \blue{smoothness and knowledge of the PDF}. The intrusive nature may lead to the loss of physical properties and suffers of the curse of dimensionality \cite{HJ, RHS}. 

\item \emph{Other methods}: \blue{moment methods} where the unknowns are the moments of the solution, \blue{stochastic collocation methods} based on  orthogonal polynomials but selecting the quadrature points \cite{NTW,ZLX}.

\end{itemize}


\subsubsection{Monte Carlo (MC) sampling methods}

Let us quickly describe the simple Monte Carlo sampling method.
Assume $u({\red z},x,t)$, $x\in\R$, solution of a PDE with uncertainty only in the initial data $u_0({\red z},x)$, $\red{ z\in\Omega\subset \R^{d_z}}$.
The method does not depend on the particular solver used for the PDE or the dimension $d_z$, and consists of three main steps. 

\begin{algorithm}[Simple Monte Carlo method]~
{\sl 
\begin{enumerate}
\item {\bf Sampling}: Sample $M$ independent identically distributed (i.i.d.) initial data 
\[
u^{k,0},\quad k=1,\ldots,M
\]
from the random initial data $u_0$ and approximate on a \blue{grid} $\Delta x$ to get $u^{k,0}_{\Delta x}$. 
\item {\bf Solving}: For each realization $u^{k,0}_{\Delta x}$ the PDE is solved by a \blue{deterministic numerical method} to obtain at time $t^n=n\Delta t$, $n>0$ 
\[
u^{k,n}_{\Delta x},\quad k=1,\ldots,M.
\]  
\item {\bf Estimating}: Estimate the desired statistical information of the quantity of interest  by its \blue{statistical mean}
\[
\E[\Psi(F(u(\cdot,t^n)))] \approx E_M[\Psi(F(u_{\Delta x}^n))]:=\frac1{M} \sum_{k=1}^M \Psi(F(u_{\Delta x}^{k,n})).
\label{mcest}
\]
\end{enumerate}
}
\end{algorithm}


In the sequel we will consider $F(u)=u$ and $\Psi(w)=w$, namely the quantity of interest is $\E[u]$.  
 Let su recall that, from the \blue{central limit theorem}, the root mean square error satisfies \cite{Caflisch, Lo77}
\be
\EE\left[(\EE[u]-E_M[u])^2\right]^{1/2} = \var(u)^{1/2}M^{-1/2}.
\ee
 Assume that the {deterministic solver} satisfies an estimate of the type
\be
\| u(\cdot,t^n)-u^n_{\Delta x}\|_{L^1(\R)} \leq C \Delta x^p,
\ee
where $p\geq 1$ and to keep notations simple we ignored the time discretization error. 
  Let us define the following \blue{norms} 
\begin{eqnarray*}
(i)&&\|\EE[u(\cdot,t^n)]-E_M[u^{n}_{\Delta x}]\|_{L^p(\R;L^2(\Omega))}:=\E\left[\left\|\EE[u(\cdot,t^n)]-E_M[u^{n}_{\Delta x}]\right\|_{L^p(\R)}^2\right]^{1/2}\\
(ii)&&\|\EE[u(\cdot,t^n)]-E_M[u^{n}_{\Delta x}]\|_{L^2(\Omega;L^p(\R))}:=\left\|\E\left[\left(\EE[u(\cdot,t^n)]-E_M[u^{n}_{\Delta x}]\right)^{2}\right]^{1/2}\right\|_{L^p(\R)}
\end{eqnarray*}
Note that, by the \emph{Jensen inequality} for any convex function $\phi$ we have 
\be
\phi\left(\E[\,\cdot\,]\right) \leq \E\left[\phi(\cdot)\right]\quad \Longrightarrow\quad (ii) \leq (i).
\ee 
Considering norm $(ii)$ and $p=1$, we have the {error estimate}
\begin{eqnarray*}
\nonumber
\|\EE[u(\cdot,t^n)]-E_M[u^{n}_{\Delta x}]\|_{L^2(\Omega;L^1(\R))}
& \leq &\underbrace{\|\EE[u(\cdot,t^n)]-E_M[u](\cdot,t^n)]\|_{L^2(\Omega;L^1(\R))}}_{\hbox{$I_1$=statistical\,\,error}}\\
&&+\underbrace{\|E_M[u(\cdot,t^n)]-E_M[u_{\Delta x}^n]\|_{L^2(\Omega;L^1(\R))}}_{\hbox{$I_2$=spatial\,\,error}}
\end{eqnarray*}
These errors are bounded by
\begin{eqnarray*}
\red{I_1} &=& \int_{\R} \EE\left[(\EE[u(\cdot,t^n)]-E_M[u(\cdot,t^n)])^2\right]^{1/2}dx \leq  C_1\sigma_u M^{-1/2} \\
\blue{I_2} &\leq & \frac1{M}\sum_{k=1}^M \int_{\R} |u^k(\cdot,t^n)-u^{k,n}_{\Delta x}| dx \leq \left(\frac1{M}\sum_{k=1}^M C_k \right)\Delta x^p = C_2\Delta x^p
\end{eqnarray*} 
with $\sigma_u=\|\var(u(\cdot,t^n))^{1/2}\|_{L^1(\R)}$. We get the final estimate
\be
\|\EE[u(\cdot,t^n)]-E_M[u^{n}_{\Delta x}]\|_{L^2(\Omega;L^1(\R))} \leq C\left({\sigma}_u M^{-1/2} + \Delta x^p\right). 
\ee
To equilibrate the discretization and the sampling errors we should take 
\[
M=O(\Delta {x}^{-2p}).
\]
Therefore, for a method of order $p$ changing the grid from $\Delta x$ to ${\Delta x}/2$ requires to multiply the number of samples by a factor $2^{2p}$.


\subsubsection{Stochastic Galerkin (SG) methods}
To describe the method, let us assume that the solution of the PDE, $u(z,x,t)$, $x\in\R$, has an uncertain initial data which depends on a one-dimensional random variable $\red{ z\in\Omega\subset \R}$. 

 The method is based on the construction of a set of orthogonal polynomials $\left\{\Phi_{m}(\theta)\right\}_{m=0}^M$, of degree less or equal to $M$, {orthonormal} with respect to the probability density function $p(z)$ \cite{Xu, PIN_book}
\[
\int_{\Omega} \Phi_n(\theta) \Phi_m(\theta) p(\theta)\,d\theta = \E[\Phi_m(\cdot)\Phi_n(\cdot)]=\delta_{mn},\qquad m,n=0,\ldots,M.
\] 
Note that, $\left\{\Phi_{m}(\theta)\right\}_{m=0}^M$ are hierarchical, in the sense that $\Phi_m(\theta)$ has degree $m$. 

The solution of the PDE is then represented as
\be
u_M(\theta,x,t)=\sum_{m=0}^M \hat{u}_m(x,t)\Phi_m(\theta),
\ee
where $\hat{u}_m$ is the {projection} of the solution with respect to $\Phi_m$
\be
\hat{u}_m(x,t)=\int_{\Omega} u(\theta,x,t)\Phi_m(\theta) p(\theta)\,d\theta=\E[u(\cdot,x,t)\Phi_m(\cdot)].
\ee 
For the quantity of interest we have
\[
\E[\Psi(F(u_M))]=\int_{\Omega} \Psi(F(u_M)) p(\theta)\,d\theta, 
\]
which can be evaluated by the same quadrature (Gaussian) used to compute $\hat{u}_m$. 

In case the quantity of interest is the \blue{expectation} of the solution we have
\[
\E[u_M]=\int_{\Omega} u_M(\theta,x,t) p(\theta)\,d\theta = \sum_{m=0}^M \hat{u}_m(x,t) \E[\Phi_m(\cdot)] =\hat{u}_0,
\]
whereas for the \blue{variance} we get
\[
\var(u_M) = \E[u_M^2]-\E[u_M]^2 = \sum_{m,n=0}^M \hat u_m \hat u_n \E[\Phi_m\Phi_n]-\hat u_0^2=\sum_{m=0}^M \hat u_m^2-\hat u_0^2.
\]
Stochastic Galerkin approximation in the field of random PDEs are better known under the name of \emph{generalized polynomial chaos (gPC)}. The solution of the PDE, is obtained by standard Galerkin approach, first replacing $u$ with $u_M$ and then projecting the PDE to the space generated by $\left\{\Phi_{m}(\theta)\right\}_{m=0}^M$. 


Let us consider a general PDE in the form
\[
\partial_t u(\theta,x,t) = {\cal J}(u(\theta,x,t)),
\]
the stochastic Galerkin method corresponds to take
\[
\E\left[\partial_t u_M(\cdot,x,t)\,\Phi_h\right] = \partial_t \hat u_h(x,t)= \E\left[{\cal J}(u_M(\cdot,x,t))\,\Phi_h\right],\quad h=0,\ldots,M.
\] 
For a {linear PDE} ${\cal J}(u(\theta,x,t))={\cal L}(u(\theta,x,t))$ we get
\[
\partial_t \hat u_h(x,t) = {\cal L}(\hat u_h (x,t)),\quad h=0,\ldots,M.
\] 
However, for nonlinear problems, for example a {bilinear PDE} $${\cal J}(u(\theta,x,t))={\cal Q}(u(\theta,x,t),u(\theta,x,t))$$ we get an additional quadratic cost $O(M^2)$ 
\[
\partial_t \hat u_h(x,t) =  \sum_{m,n=0}^M \hat u_m \hat u_n \E[{\cal Q}(\Phi_m,\Phi_n)\,\Phi_h]
,\quad h=0,\ldots,M.
\] 
A general problem, is the \blue{loss of physical properties} (like positivity of the solution or other invariants) due to the approximation in the orthogonal polynomial space.


The main interest in SG methods is due to their convergence properties, known as \emph{spectral convengence}, for smooth solutions in the random space.  
If the solution $u(\theta,x,t)\in H^r(\Omega)$, $r > 0$, the weighted Sobolev space, we have\cite{Xu}
\be
\begin{split}
\E[(u_M(\cdot,x,t)-u(\cdot,x,t))^2]^{1/2} &= \|u_M(\cdot,x,t)-u(\cdot,x,t) \|_{L^2(\Omega)}\\
& \leq C \frac{\|u(\cdot,x,t)\|_{H^r(\Omega)}}{M^r}.
\end{split}
\ee
For analytic functions, spectral convergence becomes {exponential convergence}. 

Therefore, we must equilibrate an error relation of the type
\[
M = O(\Delta x^{-p/r}),\qquad r \gg p
\]
and then very small values of $M$ are sufficient to balance the errors in the method.  

For \emph{multi-dimensional random spaces}, assuming the same degree $M$ in each dimension, the number of degrees of freedom of the polynomial space is 
\[
K=\frac{(d_z+M)!}{d_z! M!}.
\]
For example, $M=5$, $d_z=3$ gives $K=320$ (\blue{curse of dimensionality})
and typically \blue{sparse grid approximations} are necessary to avoid explosive growth in the number of parameters \cite{GS, RHS}.




\section{Uncertainty in kinetic equations}

Let us focus our attention on the specific case of kinetic equations of Boltzmann and mean-field type. More precisely, we consider kinetic equations of the general form \cite{Cer, DPR, JinPareschi, Vill}
\be
\begin{aligned}
\partial_t f({\red{z}},x,{\w},t)+ {\w}\cdot \nabla_x f({\red{z}},x,{\w},t) = \frac1{\varepsilon} Q(f,f)(\red{z},x,v,t),\quad (x,v)\in\RR^{d_x}\times \RR^{d_{\w}}
\end{aligned}
\ee
where  $\varepsilon > 0$ is the \emph{Knudsen number} and ${\red{z}\in\Omega\subseteq\RR^{d_{{z}}}}$ is a \emph{random vector}. The particular structure of the interaction term $Q(f,f)$ depends on the kinetic model considered.

Well know examples are given by the \emph{Boltzmann equation}  
\be
Q(f,f)(\red{z},x,v,t)=\int_{\mathbb{S}^{d_{\w}-1}\times\RR^{d_{\w}}}\!\!\!\!\! B({\red{z}},{\w},{\w}_*,\omega) (f(v')f(v'_*)_*-f(v) f(v_*))\,d{\w}_*\,d\omega
\label{eq:Boltzmann}
\ee
where $B({\red{z}},{\w},{\w}_*,\omega)\geq 0$ is the collision kernel and
\begin{equation}\label{eq:binary}
v^\prime = \dfrac{v+v_*}{2} + \dfrac{|v-v_*|}{2}\omega, \qquad v_*^\prime = \dfrac{v+v_*}{2} - \dfrac{|v-v_*|}{2}\omega,
\end{equation}
or by mean-field \emph{Vlasov-Fokker-Planck} type models 
\be
Q(f,f)=\nabla_{\w} \cdot \left( \mathcal P[f]f+\nabla_{\w} (D\,f) \right)
\label{eq:FP}
\ee
where $\mathcal P[\cdot]$ is a non--local operator, for example of the form 
\be
\label{eq:Bf}
\mathcal P[f]({\red{z}},x,{\w},t) = \int_{\RR^{d_x}}\int_{\RR^{d_{\w}}}P({\red{z}},x,x_*;{\w},{\w}_*)({\w}-{\w}_*)f({\red{z}},x_*,{\w}_*,t)d{\w}_*dx_*,
\ee
and $D({\red{z}},{\w},t)\ge 0$ describes the local relevance of the diffusion.



\subsection{The Boltzmann equation with random inputs}
In the classical case of rarefied gas dynamic, we have the \emph{collision invariants}
 {\be \int_{\RR^3} Q(f,f)\phi(v)\,dv
= 0,\quad \phi(v)=1,v,|v|^2,
\ee}
and in addition the entropy inequality
\be
\int_{{\RR}^3}Q(f,f)\ln(f(v))dv \leq 0. \label{eq:HTH}
\ee 
The equality holds only if $f$ is a local \emph{Maxwellian equilibrium}
\be f(\theta,v)=M(\rho,u,T)(\theta,v)=\frac{\rho(\theta)}{(2\pi T(\theta))^{d_v/2}}
\exp \left( - \frac{\vert u(\theta) - v \vert^2} {2T(\theta)}\right),
\label{eq:
MAX}
\ee
where the dependence from $x$ and $t$ has been omitted, and
\be \rho= \int_{\RR^{d_v}}f\, dv ,\quad u = {1\over
\rho} \int_{\RR^{d_v}}v f\, dv,\quad T = {1\over{d_v\rho}}
\int_{\RR^{d_v}}(v- u)^2 f\, dv, \label{eq:MQ}
\ee
are the {density}, {mean velocity} and {temperature} of the
gas depending on $({\red z},x,t)$.


Integrating the Boltzmann equation against the collision invariants $\phi(v)$ yields 
\[
\partial_t\int_{\RR^3}f({\red{z}},x,{\w},t)\phi(v)\,dv +
\nabla_x\cdot\left(\int_{\RR^3}v f({\red{z}},x,{\w},t)\phi(v)\,dv\right)
  = 0,\quad \phi(v)=1,v,|v|^2.
\label{eq:BE}
\]
These equations descrive the balance of mass, momentum and energy.
However, the system is not closed since it involves higher order moments of
{$f$}. 

The simplest way to find an approximate closure is to assume {$f\approx M$} to obtain
the compressible Euler equations with random inputs
{\begin{eqnarray}
\nonumber
&&\partial_t \rho({\red z},x,t)+{\nabla_x}\cdot (\rho u)({\red z},x,t) = 0 \\
\label{eq:euler}
\displaystyle
&&\partial_t (\rho u)({\red z},x,t)+{\nabla_x}\cdot(\rho u \otimes u +p)({\red z},x,t) = 0 \\
\displaystyle \nonumber &&\partial_t E({\red{z}},x,t)+{\nabla_x}\cdot(Eu+pu)({\red{z}},x,t) = 0,\quad p=\rho T=\frac1{d_v}(2E - 
\rho u^2).
\end{eqnarray}}
Other closure strategies, like the {Navier-Stokes} approach, lead to more accurate macroscopic approximations of the moment system.




\subsection{Numerical methods for UQ in kinetic equations}

Two peculiar aspects of kinetic equations are the \blue{high dimensionality} and the \blue{structural properties} (nonnegativity of the solution, conservation of physical quantities, $\ldots$) which represent a challenge for numerical methods. 
\red{These difficulties are even more striking in the context of UQ.} We summarize below the main advantages and drawbacks of MC and SG methods.
\medskip

\emph{MC methods for UQ} 
{
\begin{enumerate}
\item easy \blue{non intrusive} application as they rely on existing numerical solvers. Efficiency and structural properties are inherited from the existing solvers.   

\item  lower impact on \blue{the curse of dimensionality}. Easy to \blue{parallelize} and convergence is independent of the dimension of the random space.   

\item can be applied even if the PDF of the random vector is not known or lacks of regularity. 
\item {convergence behavior is slow}.  
\end{enumerate}
}

\emph{SG methods for UQ}
{
\begin{enumerate} 
\item application is \blue{intrusive} and problem dependent. Hard to combine with stochastic methods (phase space) and \red{structural properties often are lost}. 
\item suffer \blue{the curse of dimensionality}, in particular for \blue{nonlinear problems}, and special techniques are required to reduce the computational cost.  
\item require knowledge and smoothness of the PDF.
\item can achieve high accuracy, \blue{spectral accuracy} for smooth solutions, in the random space.
\end{enumerate}
In the next Sections we will focus on some of the recent progress on MC methods based on multi-fidelity techniques and on stochastic Galerkin methods using micro-macro decomposition and hybrid approaches.
}




%


\section{Single control variate (bi-fidelity) methods}

In this Section we describe the construction of multiscale control variate (MSCV) methods based on a single low fidelity model \cite{DPms}. 
 To simplify the presentation, we restrict to kinetic equations with random initial data $f(\theta,x,{\w},0)=f_0(\theta,v,{x})$ and focus on $\E[f]$ as quantity of interest. 
 
%
%
We introduce some preliminary notations. For a random variable $X$ taking values in a Banach space $\mathcal{B}({\cal R})$ we define
\[
\|X\|_{\mathcal B({\cal R}; L_2(\Omega))}=\|\E\left[X^2\right]^{1/2}\|_{\mathcal{B}({\cal R})},\qquad \|X\|_{L_2(\Omega;{\mathcal{B}({\cal R})})}=\E\left[\|X\|_{{\mathcal{B}({\cal R})}}^2\right]^{1/2}.
\] 
We assume that the equation has been discretized by a {deterministic solver} 
on a grid $\Delta {\w}$ and $\Delta x$, which satisfies \cite{DP15, RSS}
\be
\|f(\cdot,t^n)-f_{\Delta x,\Delta {\w}}^n\|_{\mathcal{B}({\cal R})} \leq C \left(\Delta x^p+\Delta {\w}^q \right),
\label{eq:det}
\ee
where, for example, $\mathcal{B}({\cal R})=L_2^1({\mathcal D}\times \R^{d_v})$ with
\be
\| f(\theta,\cdot,t)\|^p_{L^p_s(\mathcal{D}\times\RR^{d_v})} = \int_{\mathcal{D}\times\RR^{d_v}} |f(\theta,x,v,t)|^p(1+|v|^s)\,dv\,dx.
\ee 
For the Monte Carlo method therefore we have the {error estimate}
\be
\|\EE[f](\cdot,t^n)-E_M[f^{n}_{\Delta x,\Delta {\w}}]\|_{{\LB}}  
\leq C_1 \sigma_f M^{-1/2}+C_2\left(\Delta x^p + \Delta {\w}^q\right) 
\ee
with $\sigma_f = \|\var(f)^{1/2}\|_{\mathcal{B}(\mathcal{R})}$.
\subsection{Space homogeneous case}


First, we describe the method for the \emph{space homogeneous problem}
\be
\frac{\partial f}{\partial t} = Q(f,f),
\label{eq:BH}
\ee
where $f=f(\theta,\w,t)$ with initial data $f(\theta,{\w},0)=f_0(\theta,{\w})$.

Under suitable assumptions \cite{Tos1999, ToVil}, $f(\theta,{\w},t)\to f^\infty(\theta,\w)$ exponentially as as $t\to\infty$, where $f^{\infty}(\theta,{\w})$ is the \emph{Maxwellian equilibrium state} s.t. $Q(f^\infty,f^\infty)=0$. 

We denote the moments as
\[
m_{\phi}(f)(\theta,t):=\int_{\RR^{d_{\w}}}\phi({\w})f(\theta,{\w},t)d{\w},\quad \phi({\w})=1,{\w},|{\w}|^2/2.
\]
 
{\blue{MSCV methods} aim at improving the MC estimate by considering the solution of a \emph{low-fidelity model} 
$\tilde{f}(\theta,{\w},t)$, whose evaluation is significantly cheaper than $f(\theta,{\w},t)$, s.t. 
$m_\phi(\tilde{f})=m_\phi(f)$ and that
$\tilde{f}\to f^{\infty}$ as $t\to\infty$.}
The hope is that the cheaper low-fidelity model can be used
to speed up, without compromising accuracy,
the approximation of the quantities of interest corresponding to the high-fidelity model.


\subsubsection{Local equilibrium control variate}
Let us recall that a Monte Carlo estimator for $\E[f]$ based on $M$ samples gives  
\be
\|\E[f](\w,t)-E_M[f](\w,t)\|_{\LB} \leq C {\sigma}_{f} M^{-1/2}.
\ee
Now introduce the micro-macro decomposition \cite{DPZ,DPZ2, LY}  
\be
g(\theta,{\w},t)=f(\theta,{\w},t)-f^\infty(\theta,\w),
\ee
we have $g(\theta,{\w},t)\to 0$ as $t\to\infty$ and so $\var(g)\to 0$ as $t\to\infty$.

We can decompose the expected value of the solution as
\be
\begin{split}
\E[f](\w,t)=\underbrace{\E[f^\infty]({\w})}_{\rm{Accurate\,\,evaluation}}+\underbrace{\E[g]({\w},t).}_{\rm{Monte\,\,Carlo\,\,estimate}}
\end{split}
\ee
Since $f^\infty(\theta,{\w})$ is known, we can assume $\E[f^\infty]({\w})$ is evaluated with a negligible error and use the Monte Carlo estimator only on $\E[g]$ to get
\be
\begin{aligned}
&\|\E[f](\w,t)-\left(\E[f^\infty](\w)+E_M[g](\w,t)\right)\|_{\LB}\\
&\hskip 2cm =\|\E[g](\w,t)-E_M[g](\w,t)\|_{\LB} \leq C {\sigma}_{g} M^{-1/2}.
\end{aligned}
\ee
The resulting estimate now depends on ${\sigma}_{g}=\|\var(g)^{1/2}\|$ instead of ${\sigma}_{f}$,
where now ${\sigma}_g \to 0$ as $t\to\infty$. Therefore, the statistical error vanish asymptotically in time.


%
\subsubsection{Time dependent control variate}
For a time dependent low-fidelity model $\tilde{f}(\theta,v,t)$ given $M$ samples $f^{k}(v,t)$, $k=1,\ldots,M$ we define 
\be
\begin{split}
{\mathbb E}[f](\w,t) \approx {E}^{\lambda}_M[f](v,t)&:= E_M[f](v,t)-\lambda\left(E_M[\tilde f](v,t)-\E[\tilde f](v,t)\right)\\
&=\frac1{M} \sum_{k=1}^M f^{k}(v,t) - \lambda\left(\frac1{M} \sum_{k=1}^M \tilde{f}^{k}(v,t)-\tilde{\bf f}({\w},t)\right),
\end{split}
\label{eq:nest2}
\ee
with $\tilde{\bf f}=\E[\tilde f]$ or an accurate approximation. It is immediate to verify that $E^{\lambda}_M[f]$ is an unbiased estimator for any choice of $\lambda\in\RR$. 
In particular, the above estimator includes
\begin{itemize}
\item{$\lambda=0$} $\Longrightarrow$ $E^{0}_M[f]=E_M[f]$ is the {simple MC estimator}. 
\item{$\lambda=1$, $\tilde f = f^\infty$} $\Longrightarrow$  $E^{1}_M[f]$ is the {local equilibrium control variate estimator}. 
\end{itemize}
Let us consider the random variable
\be
f^{\lambda}(\theta,v,t)=f(\theta,v,t)-\lambda(\tilde{f}(\theta,v,t)-{\tilde{\bf f}}(v)).
\ee
We have $E_M^{\lambda}[f](v,t)=E_M[f^\lambda](v,t)$ and its variance is 
\[
\var(f^\lambda)=\var(f)+\lambda^2 \var(\tilde f)-2\lambda\cov(f,\tilde f).
\label{eq:var1}
\]
We can minimize the variance by direct differentiation to get
\[
\frac{\partial\var(f^\lambda)}{\partial\lambda} = 2\lambda \var(\tilde f)-2\cov(f,\tilde f)=0.
\]
As a consequence we have the following result.
\begin{proposition}
The quantity   
$\displaystyle\lambda^* = \frac{\cov(f,\tilde f)}{\var(\tilde f)}$
minimizes $\var(f^\lambda)$ at $(v,t)$ and gives
\be
\var(f^{\lambda^*}) = \left(1-\frac{\cov(f,\tilde f)^2}{\var(\tilde f)\var(f)}\right)\var(f)=(1-\rho_{f,\tilde f}^2)\var(f), 
\label{eq:var2}
\ee
where $\rho_{f,\tilde f} \in [-1,1]$ is the correlation coefficient of $f$ and $\tilde f$. We have  
\be
\lim_{t\to\infty} \lambda^*(v,t) =1,\qquad \lim_{t\to\infty} \var(f^{\lambda^*})(v,t)=0\qquad \forall\, v \in \RR^{d_v}.
\ee
\vskip -.2cm
\label{pr:1}
\end{proposition}
%
%
%
%
%
%
The resulting MSCV method can be implemented as follows.

\begin{algorithm}[space homogeneous MSCV method]~
{\sl 
\label{alg:3}
\begin{enumerate}
\item {\bf Initialize the control variate}: From the random initial data $f_0$ compute $\tilde{f}_{\Delta {\w}}^{0}$ on the mesh $\Delta {\w}$ and denote by $\tilde{\bf f}_{\Delta {\w}}^{0}$ an accurate estimate of $\E[\tilde{f}_{\Delta {\w}}^{0}]$.
\item {\bf Sampling}: Sample $M$ i.i.d. initial data $f_0^k$, $k=1,\ldots,M$ from
the random initial data $f_0$. 
\item {\bf Solving the control variate}: Compute the control variate $\tilde{f}_{\Delta {\w}}^{n}$ at time $t^n$ by a suitable scheme and denote by $\tilde{\bf f}_{\Delta {\w}}^{n}$ an accurate estimate of $\E[\tilde{f}_{\Delta {\w}}^{n}]$. 
\item {\bf Solving}: For each realization $f_0^k$, $k=1,\ldots,M$ the kinetic equation and the control variate are solved by the deterministic schemes. Denote the solutions at time $t^n$ by $f^{k,n}_{\Delta {\w}}$, and $\tilde{f}^{k,n}_{\Delta {\w}}$, $k=1,\ldots,M$. 
\item {\bf Estimating}: Estimate $\lambda^*$ using the $M$ samples as
$\displaystyle
{\lambda}_M^{*,n}= \frac{{\cov}_M(f^n_{\Delta v},\tilde{f}^{n}_{\Delta v})}{{\var}_M(\tilde{f}^{n}_{\Delta v}))}.
$
Compute the expected value of the random solution as
\vskip -.2cm
\[
{\tilde E}^{\lambda^*}_M[f^n_{\Delta v}]=\frac1{M} \sum_{k=1}^M f^{k,n}_{\Delta {\w}} - \lambda_M^{*,n}\left(\frac1{M} \sum_{k=1}^M \tilde{f}_{\Delta {\w}}^{k,n}-\tilde{\bf f}_{\Delta {\w}}^{n}\right).
\label{mcest2}
\]
\end{enumerate}
}
\end{algorithm}

Compared to standard MC no additional cost is required until the low-fidelity models can be evaluated off line, for example if we take
\[
\tilde f(\theta,\w) = \underbrace{f^\infty(\theta,\w),}_{\rm{equilibrium\,\, state}}\qquad \tilde f(\theta,\w,t) =\underbrace{ e^{-t} f_0(\theta,\w) + (1-e^{-t}) f^\infty(\theta,\w)}_{\rm{BGK\,\, approximation}}.
\]
Using such an approach one obtains an {error estimate} of the type
\begin{eqnarray}
\nonumber
\|\EE[f](\cdot,t^n)-E^{\lambda^*}_M[f^{n}_{\Delta {\w}}]\|_{{\LB}}&& \\
\nonumber
&& \hskip -2cm \leq \|\EE[f](\cdot,t^n)-E^{\lambda^*}_M[f](\cdot,t^n)\|_{{\LB}}\\
\\[-.25cm]
\nonumber
&& \hskip -1.8cm + 
\|E_M^{\lambda^*}[f](\cdot,t^n)-E^{\lambda^*}_M[f_{\Delta \w}^n]\|_{{\LB}}\\
\nonumber
&& \hskip -2cm \leq  C \left(\sigma_{f^{\lambda_*}}  M^{-1/2}+ \Delta {\w}^q\right), 
\end{eqnarray}
where $\sigma_{f^{\lambda_*}}=\|(1-\rho^2_{f,\tilde{f}})^{1/2}\var(f)^{1/2}\|_{\mathcal{B}(\mathcal{R})}$.
The statistical error depends on the correlation between $f$ and $\tilde f$. Since ${\rho}_{f,\tilde{f}}\to 1$ as $t\to\infty$ the statistical error will vanish for large times.
%


%



\begin{figure}[tb]
\begin{center}
\includegraphics[scale=0.36]{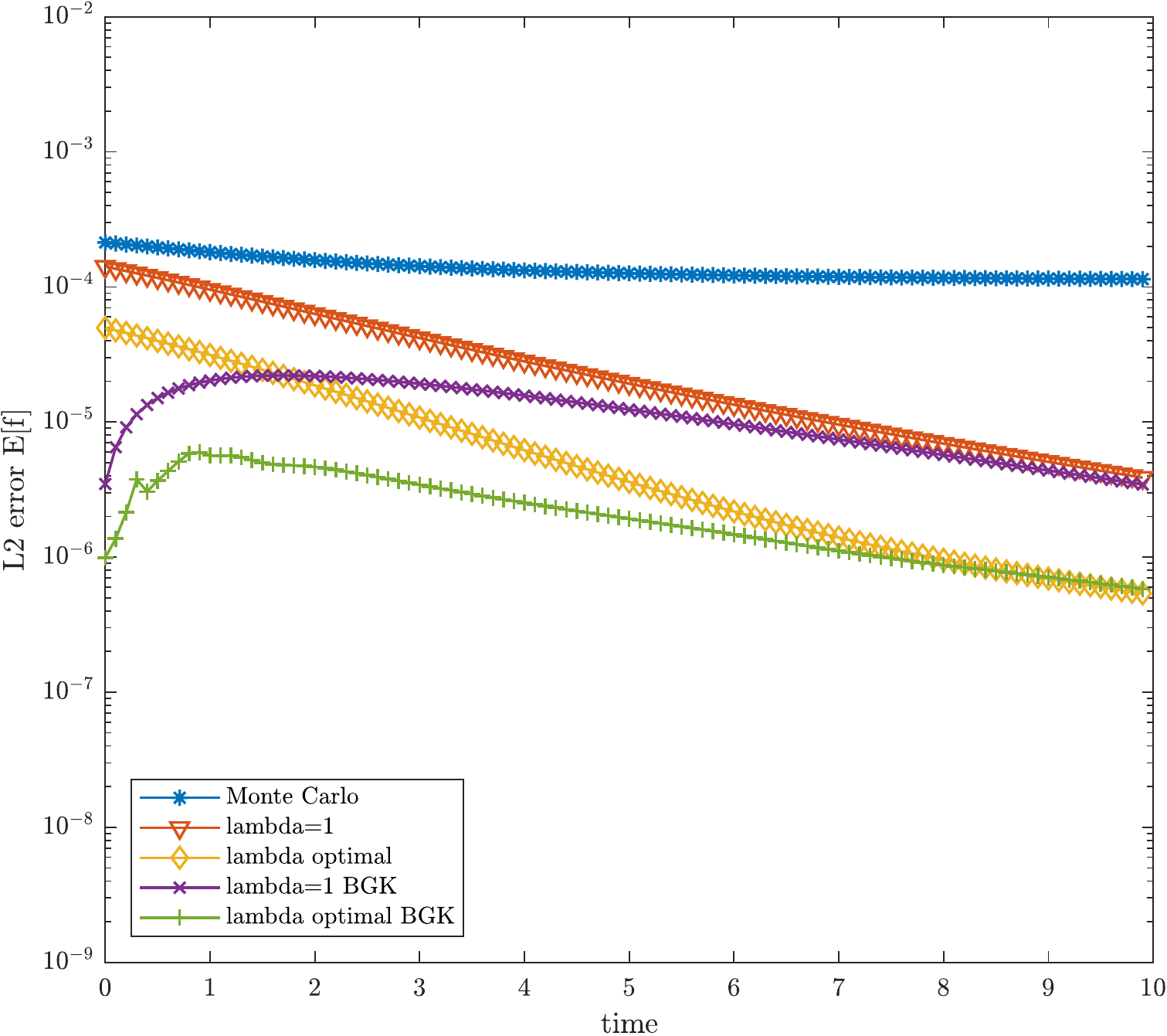}\,\,
\includegraphics[scale=0.36]{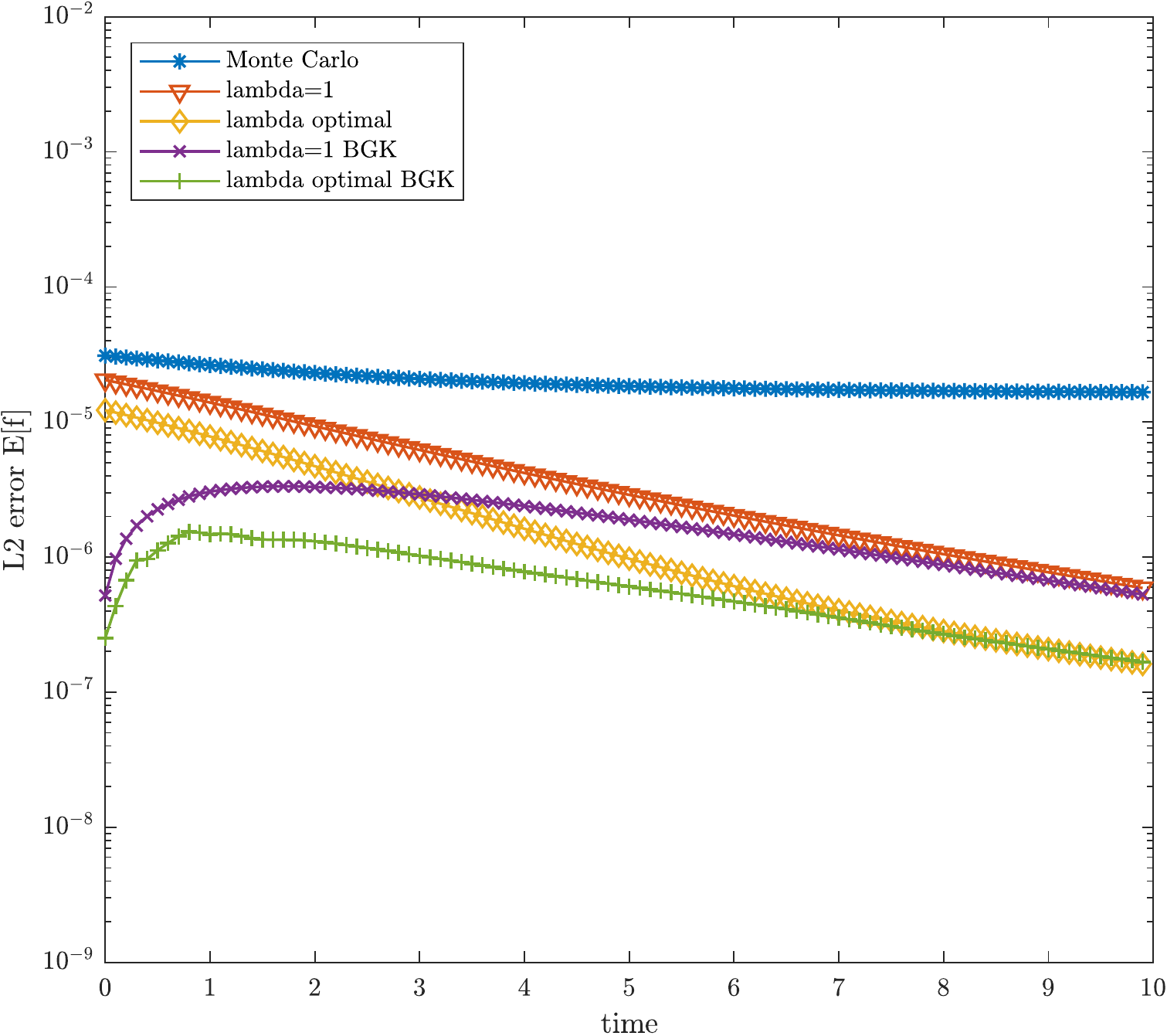}\\
\end{center}
\caption{Homogeneous relaxation. Error in $\E[f]$ using various control variate strategies. Left $M=10$, right $M=100$.}
\label{fg:uno}
\end{figure}

\begin{figure}[t]
\begin{center}
\includegraphics[scale=0.24]{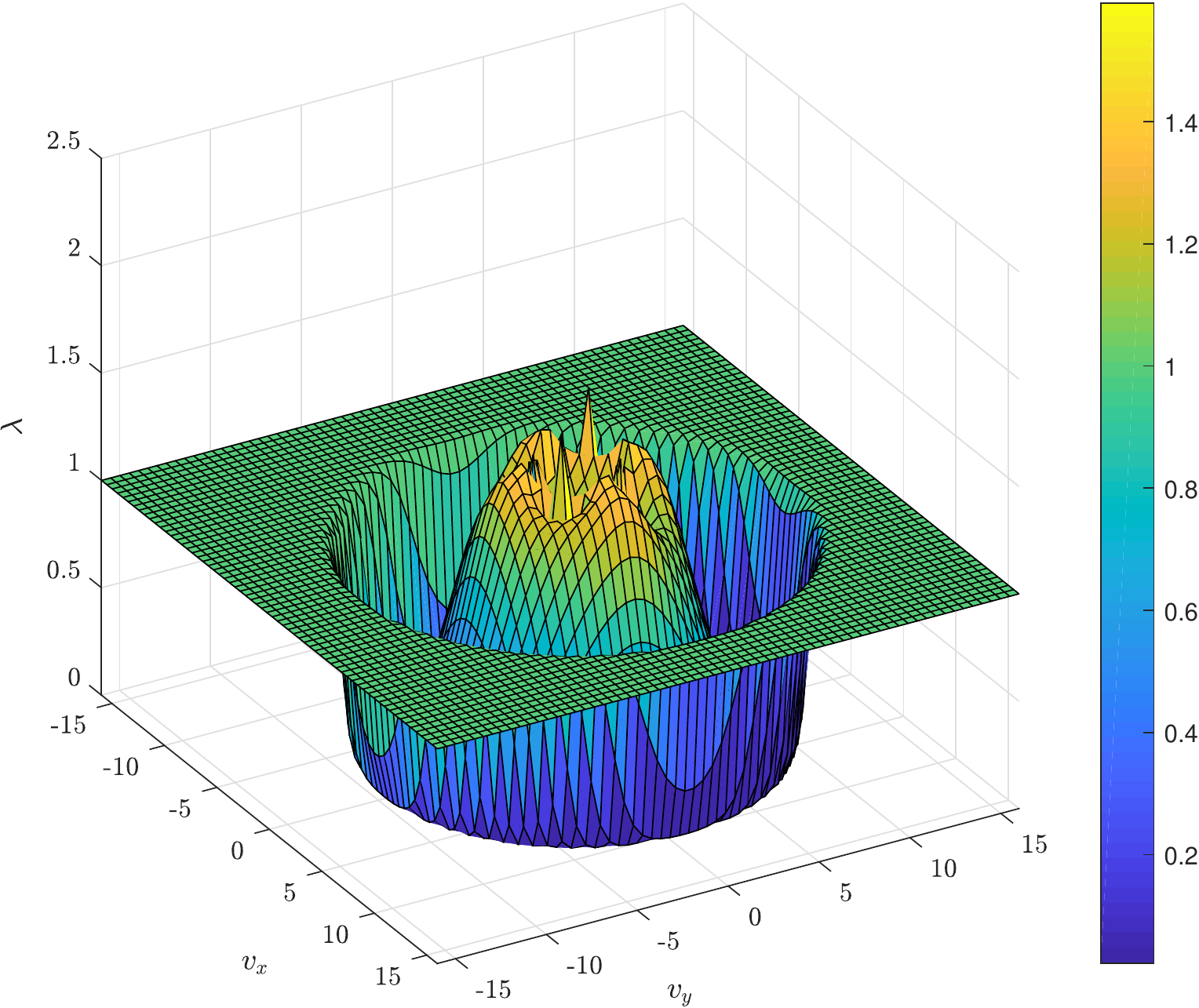}
\includegraphics[scale=0.24]{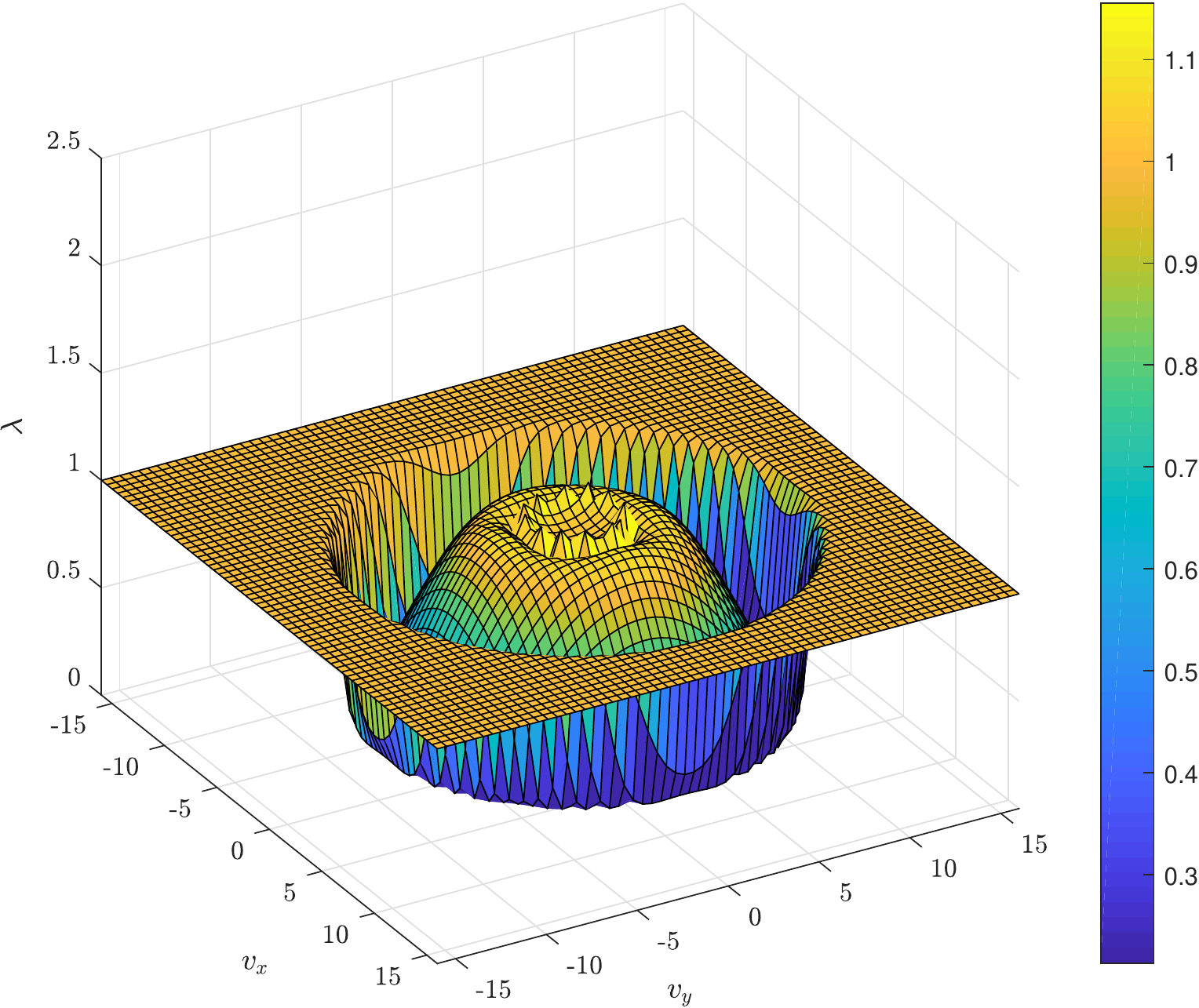}
\includegraphics[scale=0.24]{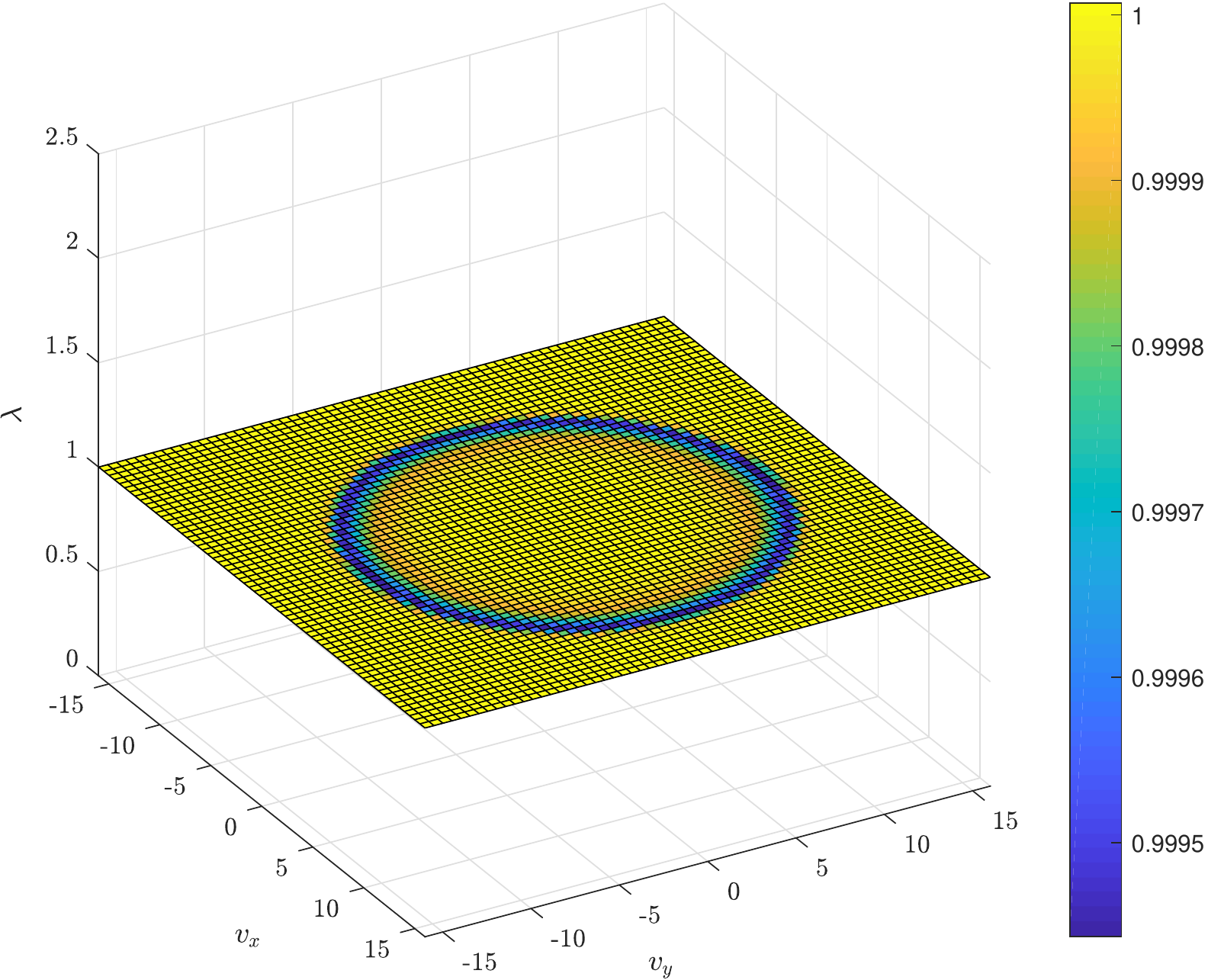}\\[-.1cm]
{\footnotesize t=5\hskip 3.3cm t=10 \hskip 3.2cm t=50}\\[+.5cm]
\end{center}
\caption{Optimal $\lambda^*(v,t)$ for the MSCV method based on a BGK control variate.}
\label{fg:due}
\end{figure}

In Figure \ref{fg:uno} we report the results for the space homogeneous Boltzmann equation with uncertain initial data using various control variates and values of $\lambda$. In Figure \ref{fg:due} we report the time evolution of the optimal value function $\lambda$ in the case of a control variate approach based on the BGK approximation.
The initial condition is a two bumps problem with uncertainty
\be
f_0(\theta,v)=\frac{\rho_0}{2\pi} \left(\exp\left(-\frac{|v-(2+s\theta)|^2}{\sigma}\right)+\exp\left({-\frac{|v+(1+s\theta)|^2}{\sigma}}\right)\right)
\ee  
with $s=0.2$, $\rho_0=0.125$, $\sigma=0.5$ and $z$ uniform in $[0,1]$. The deterministic solver adopted for the Boltzmann equation is the fast spectral method \cite{DP15, MP} and the discretization parameters are such that the stochastic error dominates the computation (see \cite{DPms} for more details). We can see from the computations that with the optimal method based on the BGK model we can gain almost two digits of precision for the same computational cost. Note that, to divide the MC error by a factor $100$ we need to multiply the number of samples by $10000$!

\subsection{Non homogeneous case}

Consider now, a general space non homogeneous kinetic equation with random inputs
\be
\partial_t f(\theta,x,v,t)+ {\w}\cdot \nabla_x f(\theta,x,v,t) = \frac1{\varepsilon} Q(f,f)(\theta,x,v,t).
\label{eq:FP_general}
\ee
For an approximated (low fidelity) solution $\tilde{f}(\theta,x,v,t)$ the estimator
reads 
\[
E^{\lambda}_M[f](x,v,t)=E_M[f](x,v,t) - \lambda\left(E_M[\tilde f](x,v,t)]-\tilde{\bf f}(x,{\w},t)\right),
\label{eq:nesth1}
\]
where $\tilde {\bf f}(x,{\w},t)$ is an accurate approximation of ${\mathbb E}[\tilde{f}](x,{\w},t)$.



The simplest control variate choice, which naturally generalizes the equilibrium control variate in the space homogeneous case, is to consider the solution of the {compressible Euler system} as control variate. Namely $\tilde{f}=f_F^\infty$ the equilibrium state corresponding to 
$U_F=(\rho_F,u_F,T_F)$ solution of the fluid model \eqref{eq:euler}. 
Improved control variates are obtained using more accurate fluid-models, like the Navier-Stokes system, or a simplified kinetic model, like a relaxation model of {BGK type} with 
\[
Q_{BGK}(\tilde{f},\tilde{f})=\nu(\tilde{f}^\infty-\tilde{f}),\quad \nu>0.
\]
The fundamental difference between the space homogeneous and the space non homogeneous case, is that now the variance of 
\be
f^{\lambda}(\theta,x,v,t)=f(\theta,x,v,t)-\lambda(\tilde{f}(\theta,x,{\w},t)-\tilde{\bf f}(x,{\w},t))
\ee
will not vanish asymptotically in time, unless the kinetic equation is close to the surrogate model (fluid regime), namely for small values of the Knudsen number.  

\begin{proposition}
The quantity   
$\displaystyle\lambda^* = \frac{\cov(f,\tilde{f})}{\var(\tilde{f})}$ 
minimizes $\var(f^\lambda)$ at $(x,v,t)$ and gives
\be
\var(f^{\lambda^*}) = (1-\rho_{f,\tilde{f}}^2)\var(f), 
\label{eq:var2h}
\ee
where $\rho_{f,\tilde{f}} \in [-1,1]$ is the correlation coefficient between $f$ and $\tilde{f}$. In addition, we have  
\be
\lim_{\varepsilon\to 0} \lambda^*(x,v,t) =1,\qquad \lim_{\varepsilon\to 0} \var(f^{\lambda^*})(x,v,t)=0\qquad \forall\, (x,v) \in \RR^{d_x}\times\RR^{d_v}.
\label{eq:alambda}
\ee
\vskip -.2cm
\label{pr:2}
\end{proposition}
Contrary to the space homogeneous case, one cannot ignore the \blue{computational cost} of solving the macroscopic fluid equations or the BGK model, although considerably smaller than that of the Boltzmann collision operator. 
Using $M_E\gg M$ samples for the control variate, we get the {error estimate} 
\bea
\nonumber
&&\|\EE[f](\cdot,t^n)-{E}^{\lambda^*}_{M,M_E}[f^{n}_{\Delta x,\Delta {\w}}]\|_{\LB}\\[-.2cm]
\label{eq:errHMMC2}
\\
\nonumber
&& \hskip 3.5cm \leq {C}\left\{\sigma_{f^{\lambda_*}} M^{-1/2}+\tau_{f^{\lambda_*}} M_E^{-1/2}+\Delta {x}^p+\Delta {\w}^q\right\} 
\eea
where $\sigma_{f^{\lambda_*}}=\| (1-\rho^2_{f,\tilde{f}})^{1/2}\var(f)^{1/2}\|_{\LL}$, $\tau_{f^{\lambda_*}}=\| \rho_{f,\tilde{f}}\var(f)^{1/2}\|_{\LL}$. 

Again the statistical error depends on the correlation between $f$ and $\tilde{f}$. In this case, ${\rho}_{f,\tilde{f}}\to 1$ as $\varepsilon\to 0$, therefore the statistical error will depend only on the fine scale sampling in the fluid limit.

We point out that, the optimal value of $\lambda$ depends on the \blue{quantity of interest} and in practice does not depend on $(x,v,t)$ unless one is interested in the details of the distribution function. For a general moment $m_\phi(f)$, the optimal value  depends on $(x,t)$ and is given by
\[
\lambda^* = \frac{\cov(m_\phi(f),m_\phi(\tilde f))}{\var(m_\phi(\tilde f))} 
\]
\begin{remark}
We remark that, by the central limit theorem we have
\[
\var(E_M[f])=M^{-1}\var(f).
\]
Therefore, taking into account the number of effective samples in the minimization process and using the independence of the estimators $E_M[\cdot]$ and $E_{M_E}[\cdot]$ we get
\bea
\nonumber
\var(E^{\lambda}_{M,M_E}[f]) &=&
M^{-1}\var(f-\lambda\tilde f)+M_E^{-1}\var(\lambda\tilde f)\\
\nonumber
&= & M^{-1}\left(\var(f)-2\lambda\cov(f,\tilde f)\right)+(M^{-1}+M_E^{-1})\lambda^2\var(\tilde f).
\eea
Minimizing with respect to $\lambda$ yields the \emph{effective optimal value} $\tilde\lambda^*$ which reads
\be
\tilde\lambda^* = \frac{M_E}{M+M_E}\lambda^*,\qquad \lambda^* = \frac{\cov(f,\tilde f)}{\var(\tilde f)}.
\ee
Let $\Co(\cdot)$ denote {computational cost} to compute the solution of a given model for a fixed value of the random parameter. The total cost is $M \Co(f) + M_E \Co({\tilde f})$. Fixing a given cost for both models $M \Co(f) = M_E \Co({\tilde f})$, we obtain 
\[
\tilde \lambda^* = \frac{\Co(f)}{\Co({f})+\Co({\tilde f})}\lambda^*.
\]
In our setting since $\Co(f) \gg \Co(\tilde f)$, or equivalently  $M_E \gg M$, we have $\tilde \lambda^* \approx \lambda^*$.
\end{remark}
\begin{figure}[tb]
\begin{center}
\includegraphics[width=5.5cm, height=4.5cm]{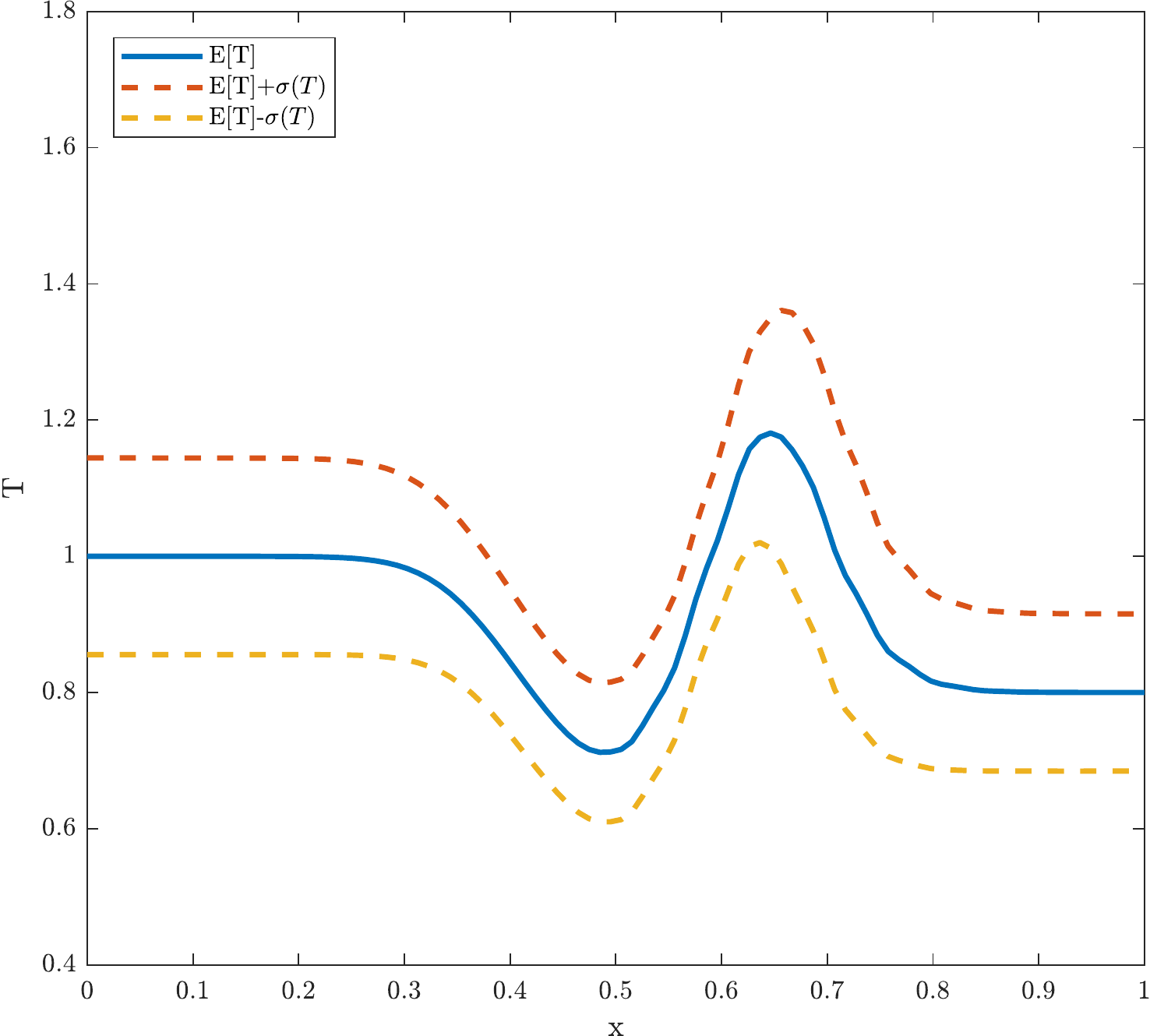}\,
\includegraphics[width=5.5cm, height=4.5cm]{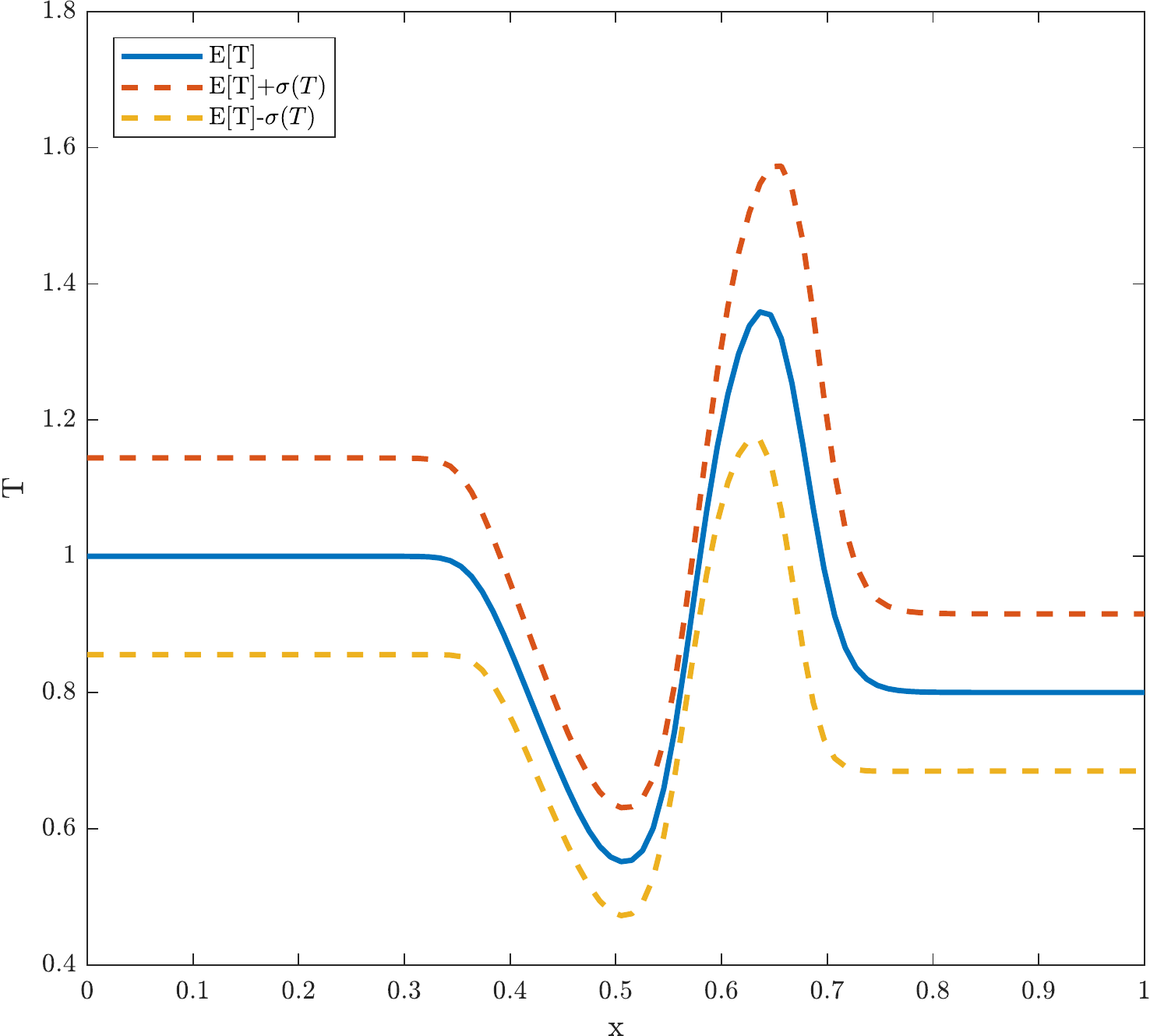}\\[+.1cm]
\includegraphics[width=5.6cm, height=4.5cm]{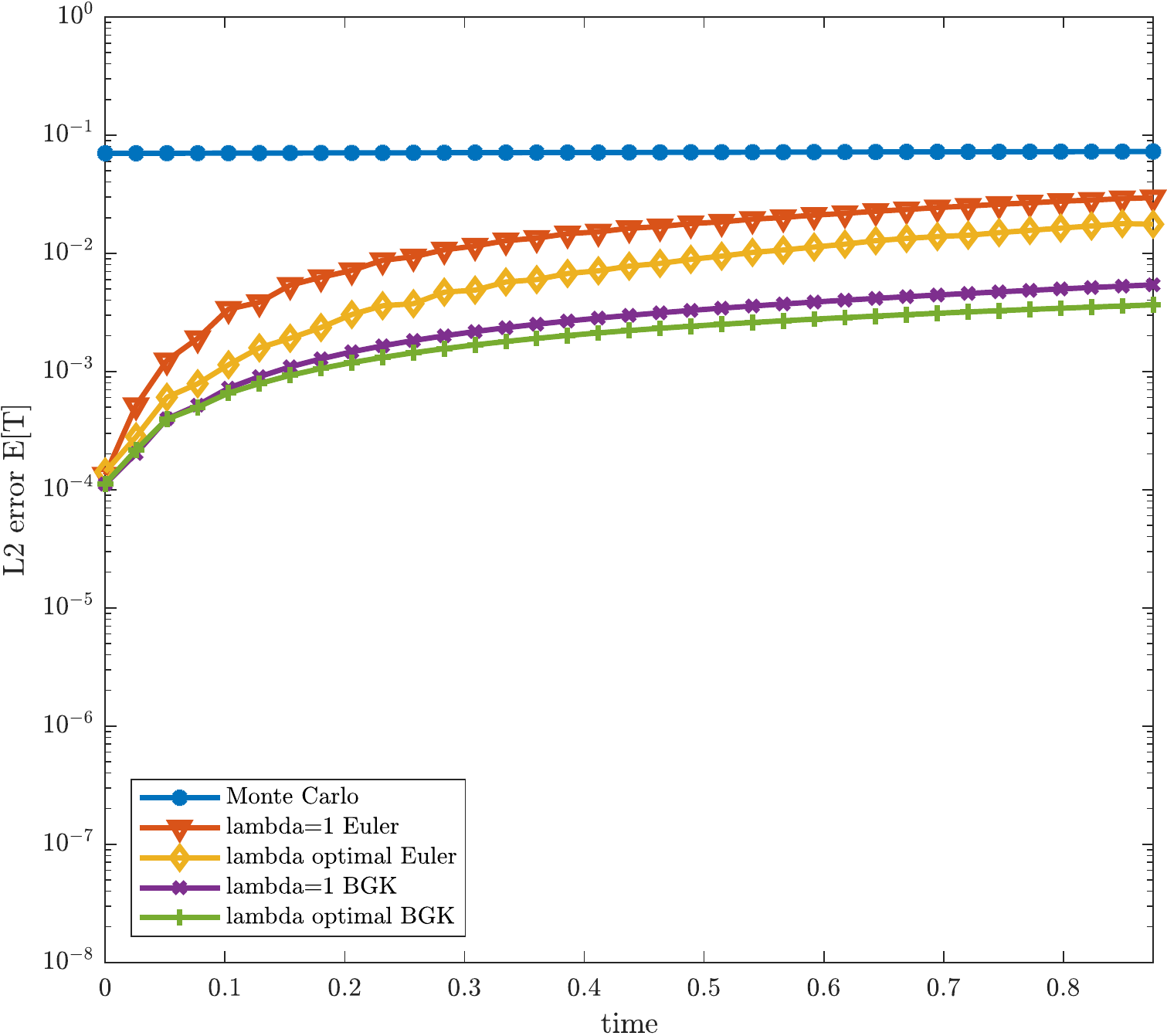}\!
\includegraphics[width=5.6cm, height=4.5cm]{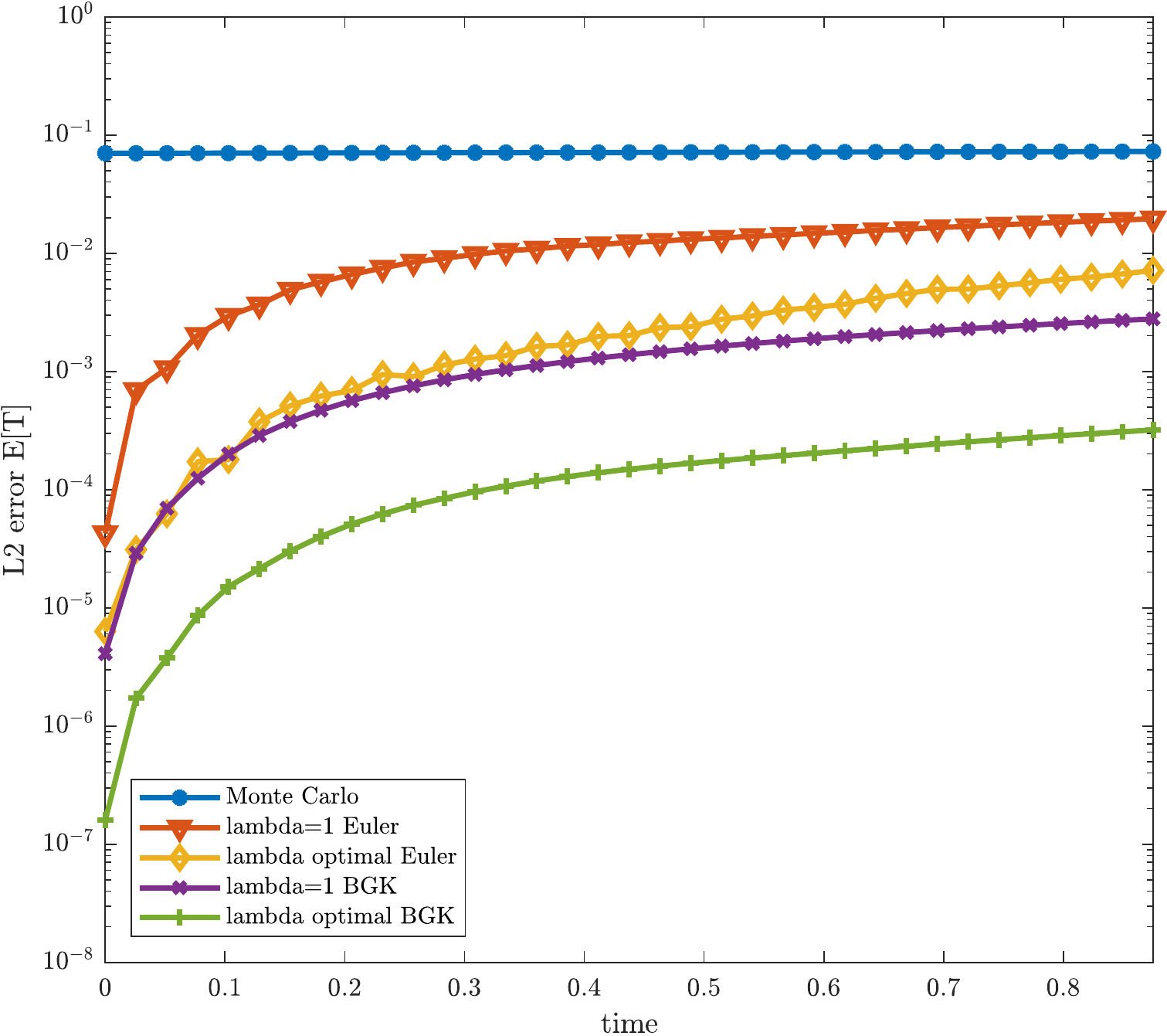}\\
\end{center}
\caption{Sod test with uncertain initial data. $\E[T]$ and confidence bands at $t=0.875$ (top). Left $\e=10^{-2}$, right $\e=10^{-3}$. Error in $\E[T]$ with $M=10$ and $\e=10^{-2}$ (bottom). Left $M_E=10^{3}$, right $M_E=10^{4}$. Here $N_x=100$, $Nv=32$.}
\label{fg:tre}
\end{figure}
As a numerical example of the performance of the method let us consider the Boltzmann equation with $d_x=1$, $d_v=2$ for the following Sod test with uncertain initial data
\be
\begin{split}
&\rho_0(x)=1, \ \ T_0(\theta,x)=1+s\theta \qquad &\textnormal{if} \ \ 0<x<L/2 \\
&\rho_0(x)=0.125,\ T_0(\theta,x)=0.8+s\theta  \qquad &\textnormal{if} \ \ L/2<x<1
\end{split}
\label{eq:sod}
\ee
with $s=0.25$, $\theta$ uniform in $[0,1]$ and equilibrium initial distribution
\[
f_0(\theta,x,\w)=\frac{\rho_0(x)}{2\pi } \exp\left({-\frac{|v|^2}{2T_0(\theta,x)}}\right).
\]
Since we are interested only in the accuracy in the random variable, the numerical parameters of the deterministic discretization $N_v$, $N_x$ and $\Delta t$ have been selected such that the deterministic error is smaller than the stochastic one (see \cite{DPms} for further details). In Figure \ref{fg:tre} (top), we report the expectation of the solution at the final time together with the confidence bands. In the same Figure (bottom) we also report the various errors using different control variates for the expected value of the temperature as a function of time. The optimal values of $\lambda^*(x,t)$ have been computed with respect to the temperature. The improvements obtained by the various control variates are evident and, as expected, becomes particularly striking close to fluid regimes.
\begin{figure}[tb]
\begin{center}
\includegraphics[width=5.5cm, height=4.5cm]{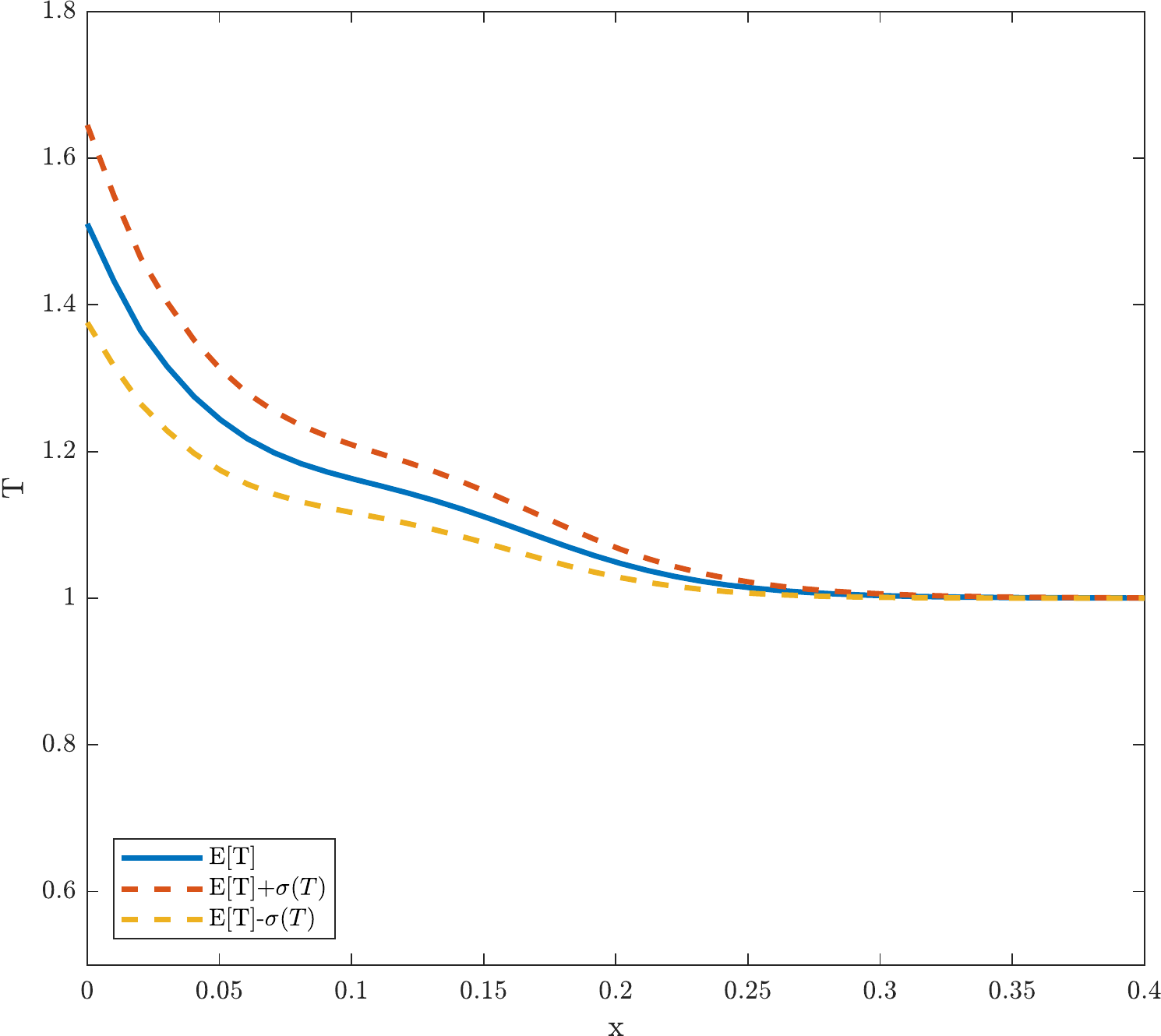}\,
\includegraphics[width=5.5cm, height=4.5cm]{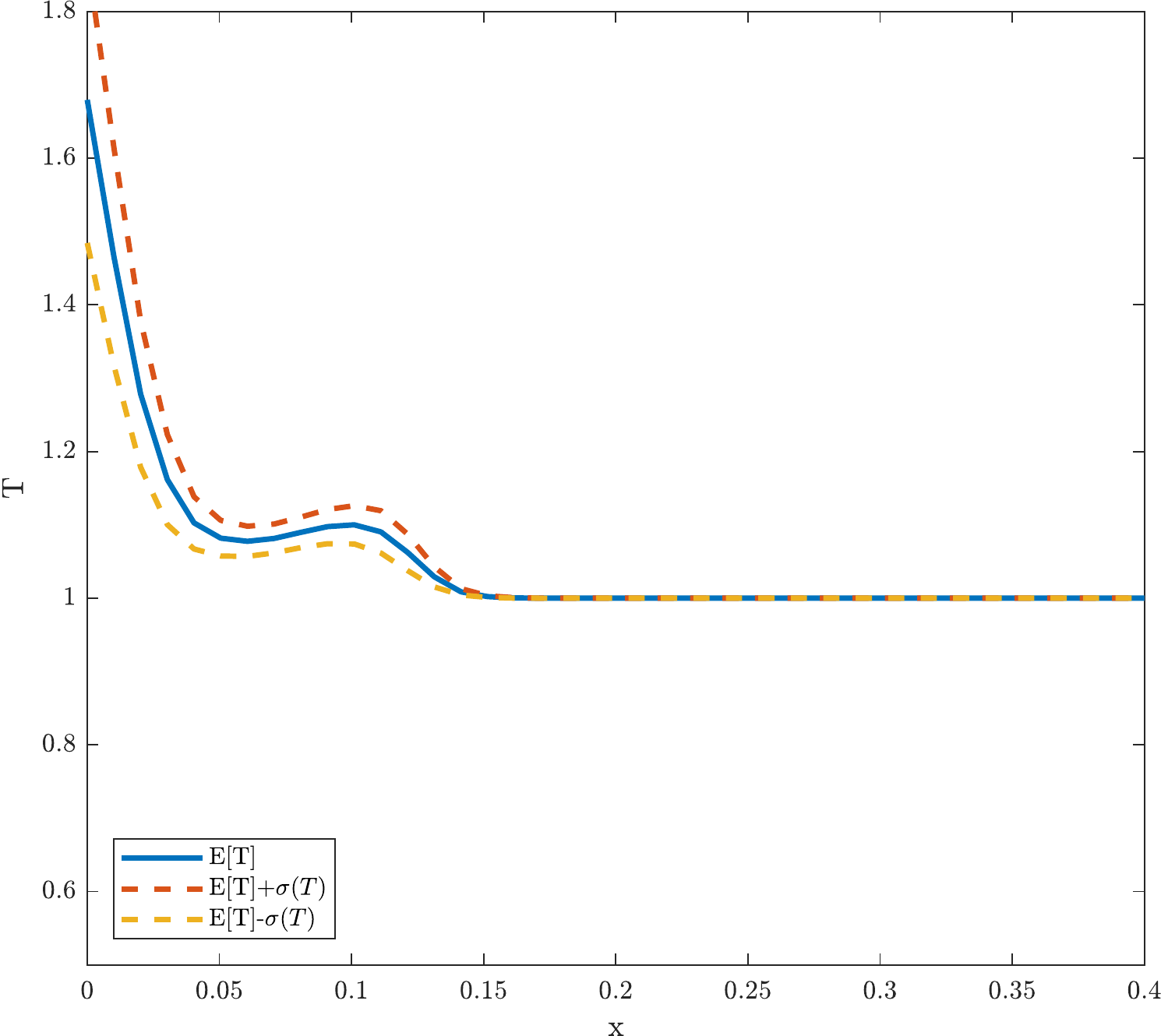}\\[+.1cm]
\includegraphics[width=5.6cm, height=4.5cm]{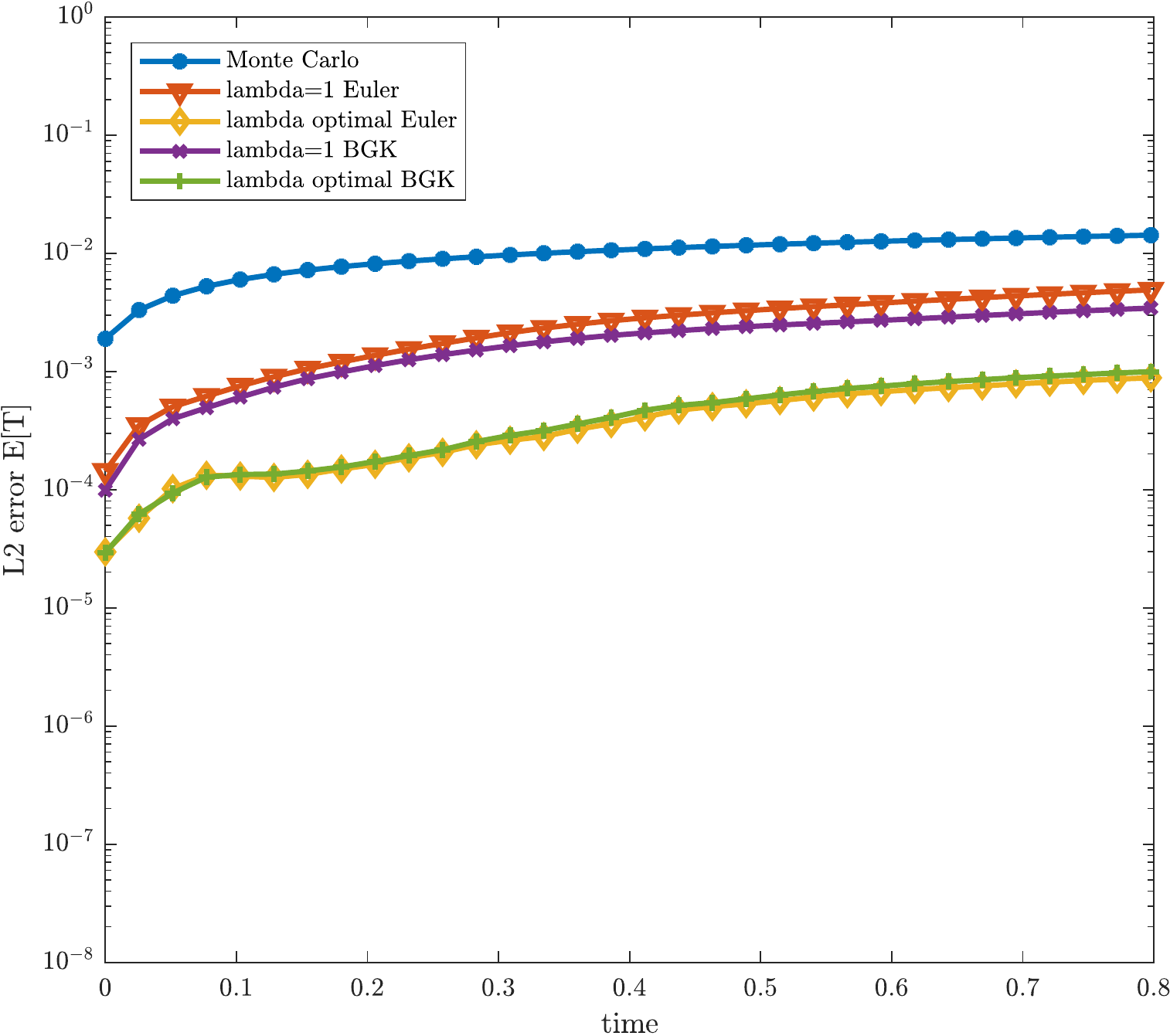}\,
\includegraphics[width=5.6cm, height=4.5cm]{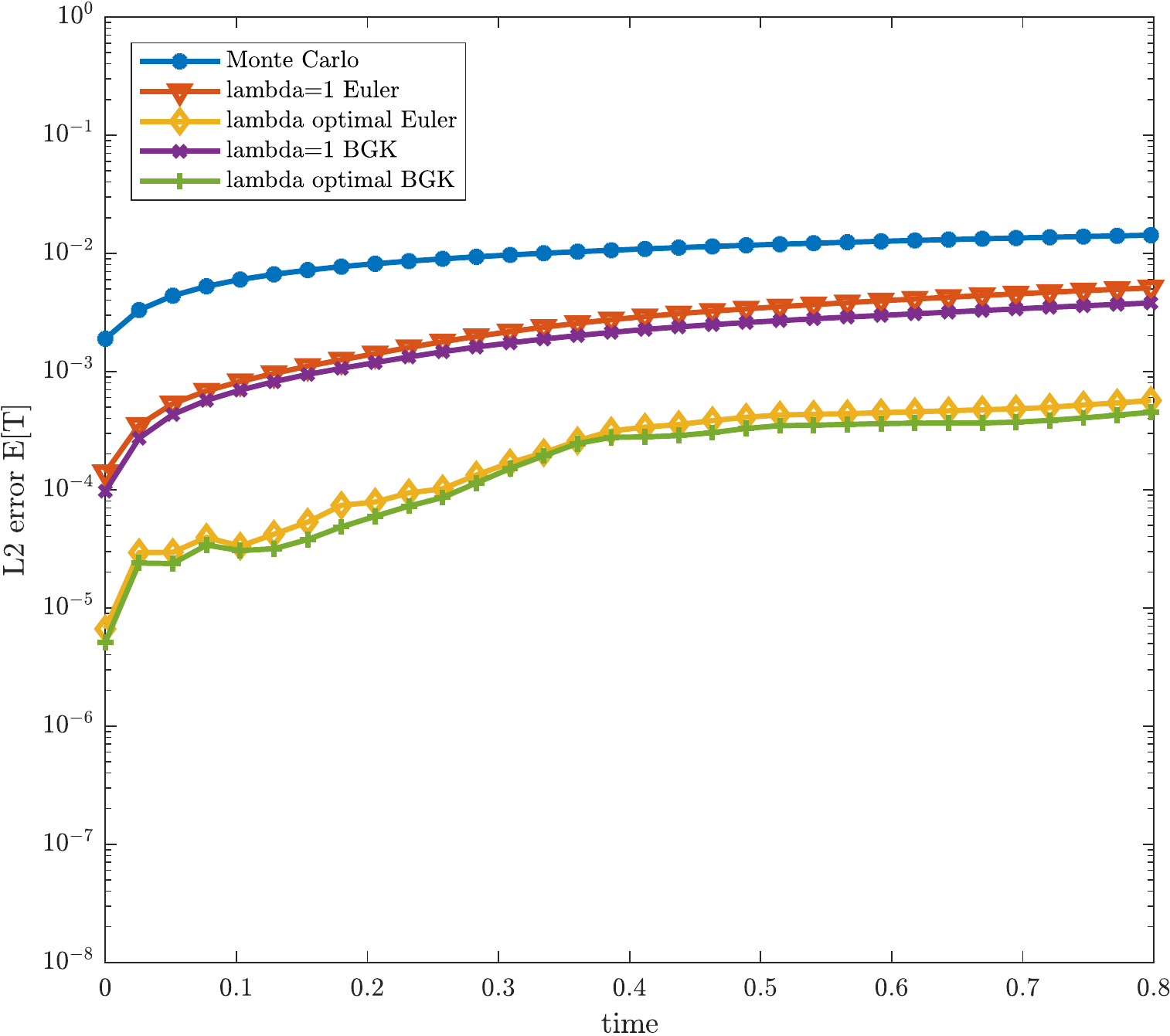}\\
\end{center}
\caption{Sudden heating with uncertain boundary condition. $\E[T]$ and confidence bands at $t=0.9$ (top). Left $\e=10^{-2}$, right $\e=10^{-3}$. Error in $\E[T]$ with $M=10$ and $\e=10^{-3}$ (bottom). Left $M_E=10^{3}$, right $M_E=10^{4}$. Here $N_x=100$, $Nv=32$.}
\label{fg:quattro}
\end{figure}

Next we consider in the same setting the sudden heating problem with uncertain boundary condition.
Initial condition is a local equilibrium with $\rho_0=1$, $u_0=0$, $T_0=1$ for $x\in[0,1]$ with diffusive equilibrium boundary condition with uncertain \blue{wall temperature}
\bea
T_w(\theta,0)=2(T_0+s\theta)
\label{eq:shp}
\eea
with $s=0.2$, $\theta$ uniform in $[0,1]$.
The results are summarized in Figure \ref{fg:quattro}, where we report the expectation of the temperature at the final time and the various errors using different control variates. In this case, due to the source of uncertainty at the boundary there is no relevant difference between the Euler and BGK control variates and the results is less sensitive to the choice of the Knudsen number.



\section{Multiple control variate (multi-fidelity) methods}
\begin{figure}[tb]

\begin{center}
\includegraphics[scale=0.25]{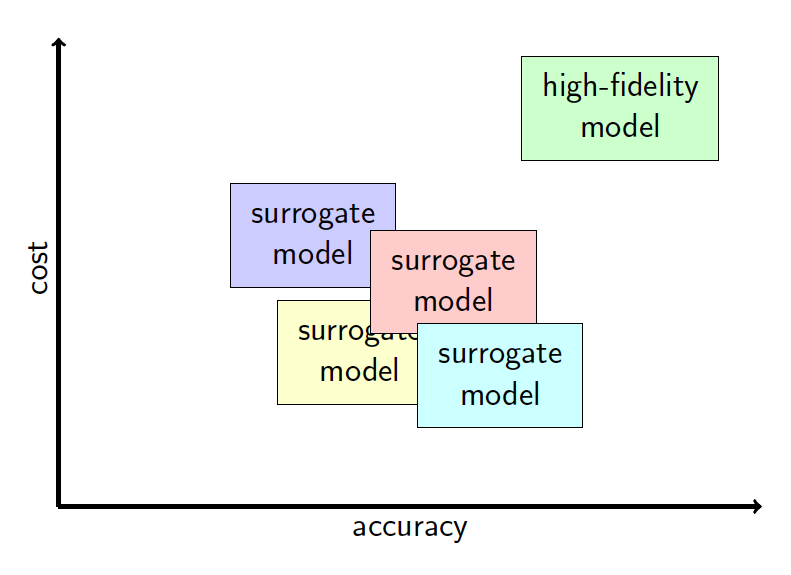}
\end{center}
\caption{High fidelity model and low fidelity surrogate models vs cost}
\label{fg:mcv}
\end{figure}

The bi-fidelity approach developed in the previous Sections is fully general and accordingly to the particular kinetic model studied one can select a suitable approximated solution as control variate which acts at a given scale. In this section we extend the methodology to the use of several approximated solutions as control variates with the aim to further improve the variance reduction properties of MSCV methods (see Figure \ref{fg:mcv}).


Given $\tilde{f}_{1},\ldots,\tilde{f}_{L}$ approximations of $f(\theta,v,t)$ we can consider the random variable
\be
f^{\lambda_1,\ldots,\lambda_L}(\theta,v,t)=f(\theta,v,t)-\sum_{h=1}^L \lambda_h (\tilde{f}_{h}(\theta,v,t)-\tilde{\bf f}_{h}(v,t)),
\ee  
where $\tilde{\bf f}_{h}(v,t)=\E[\tilde{f}_{h}](v,t)$. 
 
The variance is given by
\begin{eqnarray*}
\var(f^{\lambda_1,\ldots,\lambda_L})&=&\var(f)+\sum_{h=1}^L \lambda_h^2 \var(\tilde{f}_{h})\\
&&+2\sum_{h=1}^L\lambda_h\left(\sum_{k=1 \atop k\neq h}^L \lambda_k
\cov(\tilde{f}_{h},\tilde{f}_{k}) - \cov(f,\tilde{f}_{h})\right),
\end{eqnarray*} 
\vskip -.5cm
or in vector form
\be
\var(f^{\llambda})=\var(f)+\llambda^T  C \llambda - 2\llambda^T b
\label{eq:varmcv}
\ee
where ${\llambda}=(\lambda_1,\ldots,\lambda_L)^T$, ${b}=(\cov(f,{f}_{1}),\ldots,\cov(f,{f}_{L}))^T$ and
$C=( c_{ij})$, ${ c_{ij}}=\cov({f}_{i},{f}_{j})$ is the symmetric $L\times L$ \emph{covariance matrix}.


\begin{proposition}
Assuming the covariance matrix $C$ is not singular, the vector
\be
\llambda^* = {C}^{-1} b,
\label{eq:lambdamcv}
\ee
minimizes the variance of $f^\llambda$ at the point $(v,t)$ and gives
\be
\var(f^{\llambda^*})= \left(1-\frac{b^T(C^{-1})^T b}{\var(f)}\right)\var(f).
\label{eq:varmvcs}
\ee
\label{th:2}
\end{proposition}

In fact, the optimal values $\lambda^*_h$, $h=1,\ldots,L$ are found by equating to zero the partial derivatives with respect to $\lambda_h$. 
%
%
This corresponds to the linear system
\be
\cov(f,\tilde{f}_{h}) = \sum_{k=1}^L \lambda_k \cov(\tilde{f}_{h},\tilde{f}_{k}),\quad h=1,\ldots, L,
\ee
or equivalently $C\Lambda = b$.

\begin{figure}[tb]
\begin{center}
\includegraphics[scale=0.5]{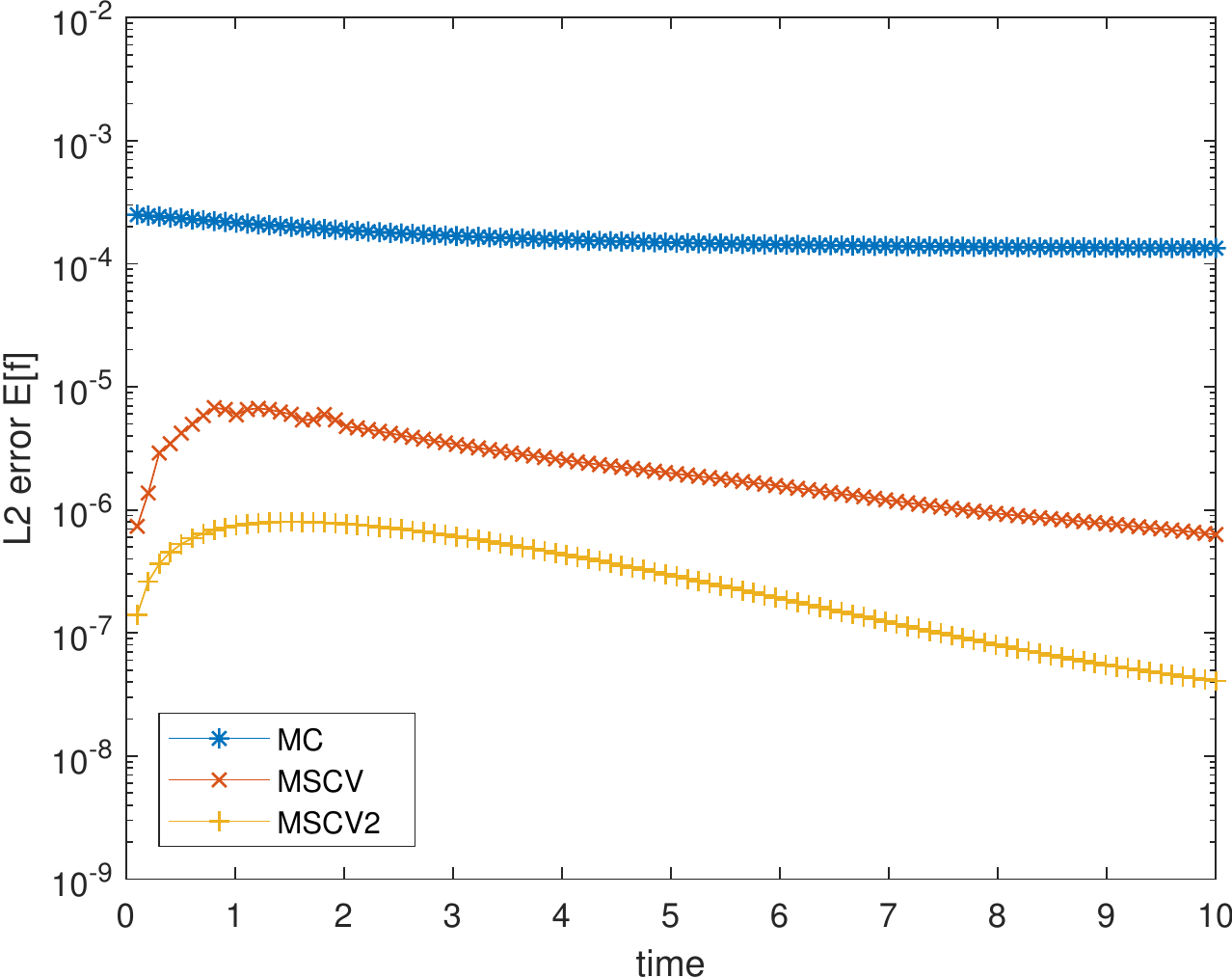}
\end{center}
\caption{Homogeneous relaxation. Error for $\E[f]$ over time with $M=10$ for the MSCV method based on BGK and the MSCV2 method based on the two control variates $f_0$ and $f^\infty$.}
\label{fg:hom2}
\end{figure}

\begin{example}
Let us consider the case $L=2$, where $\tilde f_1 = f_0$ and $\tilde f_2 = f^\infty$. \\

The optimal values $\lambda^*_1$ and $\lambda_2^*$ are readily found and are given by
\begin{eqnarray}
\nonumber
\lambda_1^* &=& \frac{\var(f^\infty)\cov(f,f_0)-\cov(f_0,f^\infty)\cov(f,f^\infty)}{\Delta},\\[-.1cm]
\label{eq:l12}
\\[-.2cm]
\nonumber
\lambda_2^* &=& \frac{\var(f_0)\cov(f,f^\infty)-\cov(f_0,f^\infty)\cov(f,f_0)}{\Delta},
\end{eqnarray}
where $\Delta = \var(f_0)\var(f^\infty)-\cov(f_0,f^\infty)^2$. 

Using $M$ samples the optimal estimator reads
\be
\begin{split}
E^{\lambda^*_1,\lambda^*_2}_M(v,t) =&E_M[f](v,t)-\lambda^*_1 \left(E_M[f_0](v)-{\bf f}_0(v)\right)\\
&-\lambda^*_2 \left(E_M[f^\infty](v)-{\bf f}^\infty(v)\right).
\end{split}
\label{eq:nest2b}
\ee
 
Since $\displaystyle\lim_{t\to\infty }f(v,t)= f^\infty(v)$ we get  
\[
\lim_{t\to\infty} \lambda_1^* = 0, \qquad \lim_{t\to\infty} \lambda_2^* = 1,
\]
and thus, the variance of the estimator vanishes asymptotically in time
\[
\lim_{t\to\infty} E^{\lambda^*_1,\lambda^*_2}_M(v,t) = {\bf f}^\infty(v).
\]
\end{example}
In Figure \ref{fg:hom2} we report the results obtained for the homogeneous relaxation problem with uncertain initial data. Compared to the optimal BGK control variate, at the same computational cost, using the estimator based on two control variates described above we can gain one additional digit of accuracy. 


%

\subsection{Hierarchical methods}

Now, let us assume $f_1,\ldots,f_L$ control variates with an \blue{increasing level of fidelity}. The idea is to apply recursively a bi-fidelity approach where the level $f_{h-1}$ is used as control variate for the level $f_{h}$.

To start with, we estimate $\EE[f]$ with $M_L$ samples using $f_L$ as control variate
\[
\EE[f]\approx E_{M_L}[f]-{\hat \lambda}_L\left(E_{M_L}[f_L]-\EE[f_L]\right).
\] 
Next, to estimate $\EE[f_L]$ we use $M_{L-1}\gg M_L$ samples with $f_{L-1}$ as control variate 
\[
\EE[f_L]\approx E_{M_{L-1}}[f_L]-{\hat \lambda}_{L-1}\left(E_{M_{L-1}}[f_{L-1}]-\EE[f_{L-1}]\right).
\]
 
Similarly, in a \blue{recursive way}, we can construct estimators for the remaining expectations of the control variates $\EE[f_{L-1}],\EE[f_{L-2}],\ldots,\EE[f_2]$ using respectively $M_{L-2}\ll M_{L-3}\ll\ldots \ll M_1$ samples
until 
\[
\EE[f_2]\approx E_{M_{1}}[f_2]-{\hat \lambda}_{1}\left(E_{M_{1}}[f_{1}]-\EE[f_{1}]\right),
\]
and we stop with the final estimate using $M_0 \gg M_1$
\[
\EE[f_1]\approx E_{M_0}[f_1].
\]


By combining the estimators of each stage we define the \emph{hierarchical estimator} 
\begin{eqnarray}
\nonumber
E_L^{{\hat \llambda}}[f] &:=& E_{M_L}[f] -{\hat \lambda}_L\left(E_{M_L}[f_L]-E_{M_{L-1}}[f_L]\right.\\
\nonumber
&+&{\hat \lambda}_{L-1}\left(E_{M_{L-1}}[f_{L-1}]-E_{M_{L-2}}[f_{L-1}]\right.\\[-.3cm]
\label{eq:lambdar}
\\[-.2cm]
\nonumber
&\ldots&\\
\nonumber
&+&{\hat \lambda}_{1}\left.\left.\left(E_{M_{1}}[f_{1}]-E_{M_{0}}[f_{1}]\right) \ldots\right)\right).
\end{eqnarray} 
The estimator can be recast in the form 
\bea
\nonumber
E_L^{{\hat \llambda}}[f] &=&  E_{M_{L}}[f_{L+1}]-\sum_{h=1}^{L}\lambda_h(E_{M_h}[f_h]-E_{M_{h-1}}[f_h])\\[-.2cm]
\label{eq:mscvgr}
\\[-.1cm]
\nonumber
&=&\lambda_1 E_{M_0}[f_1] + \sum_{h=1}^{L} (\lambda_{h+1} E_{M_h}[f_{h+1}]-\lambda_{h} E_{M_{h}}[f_{h}]),
\eea
where we defined $\hat\llambda =(\hat \lambda_1,\ldots,\hat \lambda_L)^T$ and
\be
\lambda_h=\prod_{j=h}^L \hat\lambda_j,\quad h=1,\ldots,L,\qquad \lambda_{L+1}=1,\qquad f_{L+1}=f.
\label{eq:lambdah}
\ee


The total variance of the resulting estimator is 
\bea
\nonumber
&&\var(E_L^{{\hat \llambda}}[f])=\lambda_1^2 M_0^{-1} \var(f_1) \\[-.2cm]
\\[-.2cm]
\nonumber
&&\hskip 1cm + \sum_{h=1}^{L} M_h^{-1}\left\{\lambda_{h+1}^2\var(f_{h+1})+\lambda_{h}^2\var(f_{h})-2\lambda_{h+1}\lambda_{h}\cov(f_{h+1},f_h)\right\}.
\eea 
By direct differentiation we get the \blue{tridiagonal system} for $h=1,\ldots,L$
\be
\begin{split}
&M_{h-1}^{-1}\left\{\lambda_h\var(f_h)-\lambda_{h-1} \cov(f_h,f_{h-1})\right\}\\
&+ M_{h}^{-1}\left\{\lambda_h \var(f_h)-\lambda_{h+1}\cov(f_{h+1},f_{h})\right\}=0,
\end{split}
\ee
which under the assumption $M_h \ll M_{h-1}$ leads to the \emph{quasi-optimal} solutions
\be
\lambda^*_h=\prod_{j=h}^L \hat\lambda^*_j, \qquad {\hat \lambda}_j^* = \frac{\cov(f_{j+1},f_{j})}{\var(f_{j})}.
\ee


In the case of a space homogeneous kinetic equation the hierarchical MSCV estimator \eqref{eq:mscvgr} satisfies the {error bound} 
\bea
\nonumber
&&\|\EE[f](\cdot,t^n)-{E}^{\hat\Lambda^*}_{L}[f^n_{\Delta v}]\|_{{\LB}} \\[-.2cm]
\label{eq:errREC2}
\\[-.2cm]
\nonumber
&& \hskip 4cm \leq C \left(\sum_{h=1}^L \xi_h \sigma_h M_h^{-1/2}+\xi_0 M_0^{-1/2}+\Delta v^q \right)
\eea
where $\sigma_h=\left\|\left(1-\rho^2_{f_h,f_{h-1}}\right)^{1/2}\var(f_h)^{1/2}\right\|_{\LL}$, $\tau_h=\|\rho_{f_h,f_{h-1}}\var(f_h)^{1/2}\|_{\LL}$ and $\xi_h = \prod_{j=h+1}^L \tau_j$.

 If the control variates share the same behavior as $t\to \infty$, namely $f_h\to f^{\infty}$ for $h=1,\ldots,L$, we get $\rho^2_{f_h,f_{h-1}}\to 1$ as $t\to \infty$ the statistical error depends only on the finest level of samples $M_0$. Similar considerations hold in the space non homogeneous case as $\e\to0$.
 
\begin{figure}[tb]
\begin{center}
\includegraphics[scale=0.4]{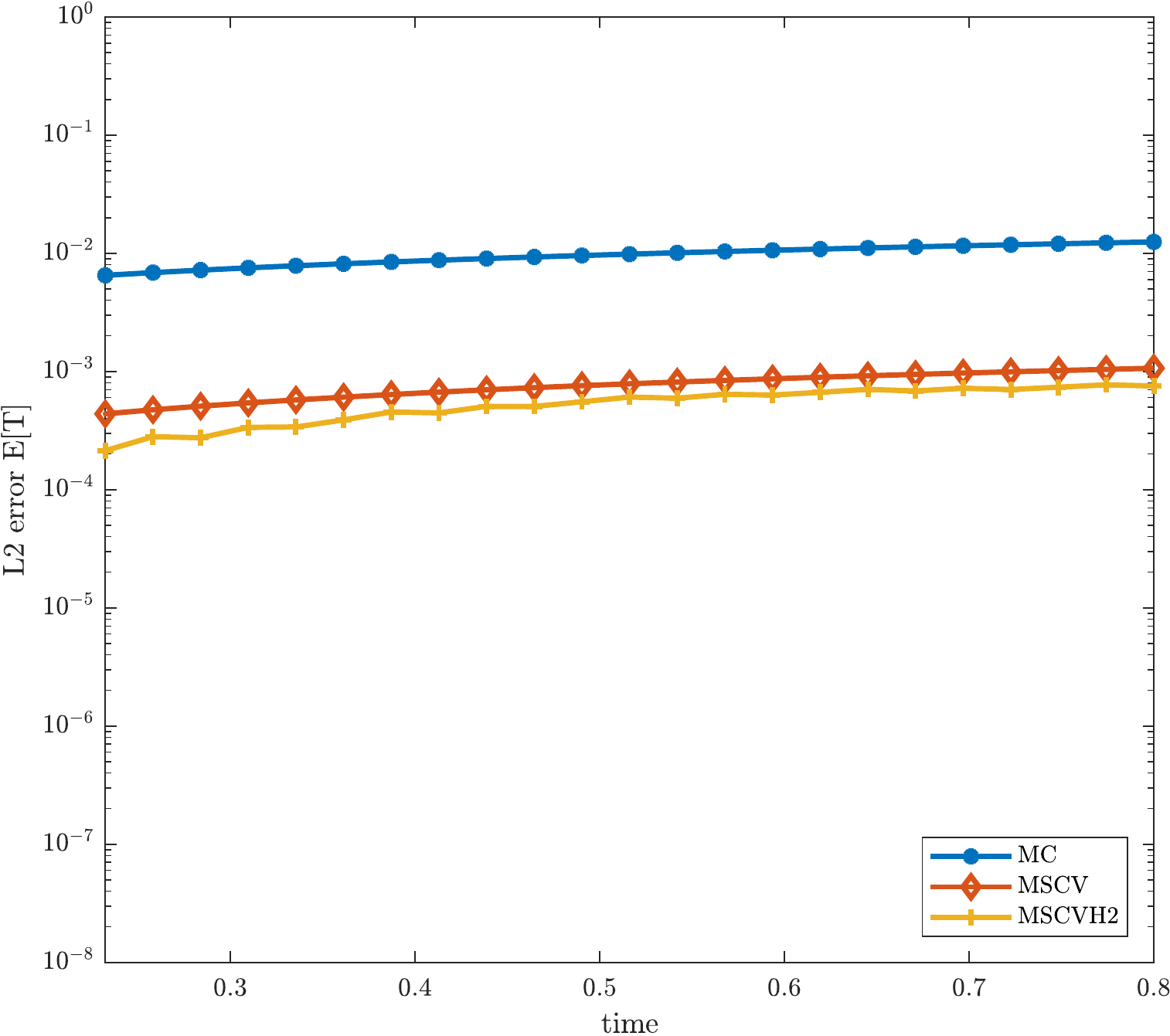}
\includegraphics[scale=0.4]{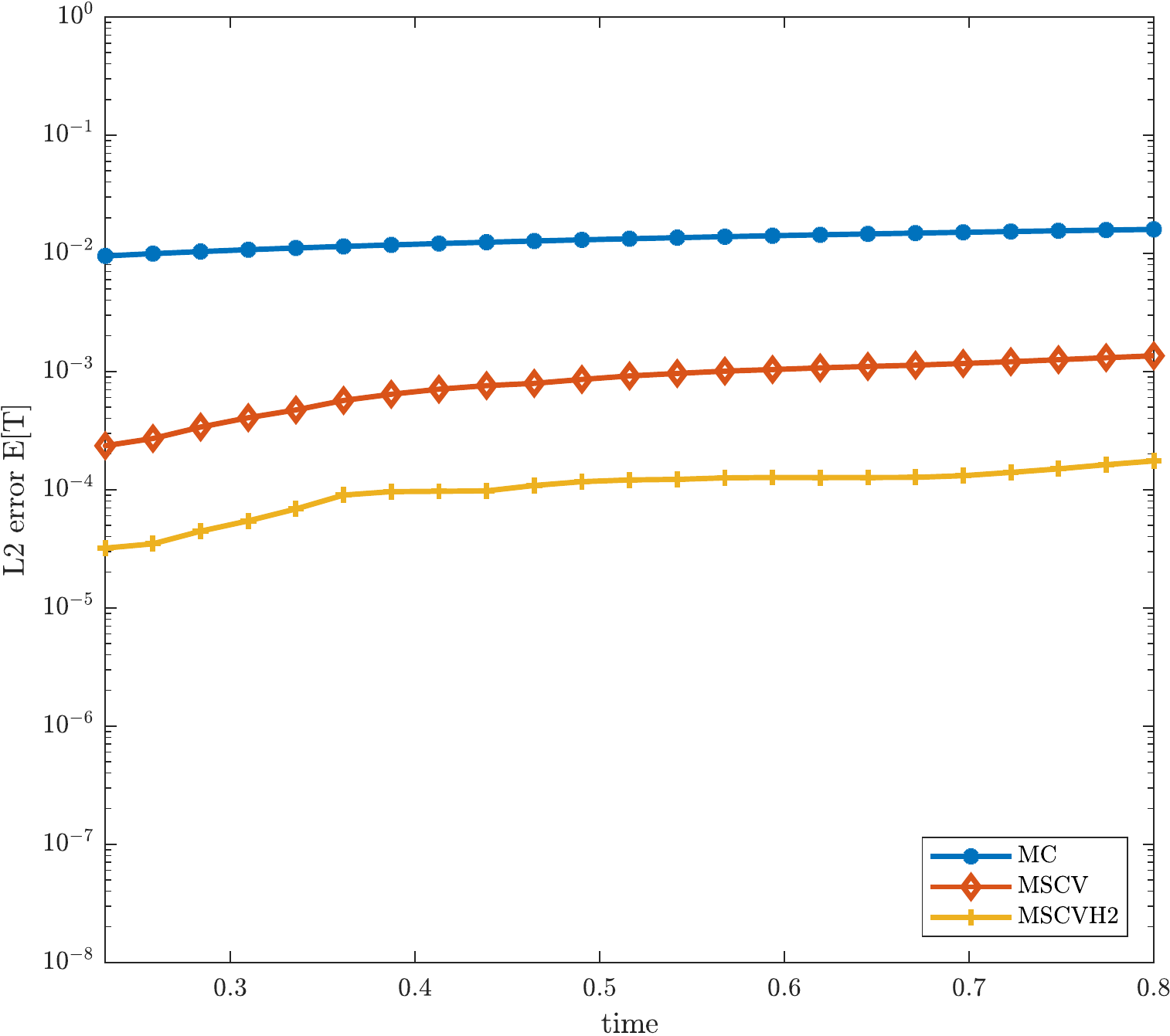}
\end{center}
\caption{Sudden heating with uncertain boundary condition. {Error for $\E[T]$ over time for $\varepsilon=10^{-2}$ (left) and $\varepsilon=10^{-3}$ (right). MSCV method based on BGK and MSCVH2 based on BGK and Euler. $M_2=10$ for Boltzmann, $M_1=10^2$ for BGK and $M_0=10^4$ for Euler.}}
\label{fig1Mc}
\end{figure}

In Figure \ref{fig1Mc} the results obtained in the case of the sudden heating problem \eqref{eq:shp} using a three models hierarchy based on the Euler system, the BGK model and the full Boltzmann equation are reported. 



\subsection{Multi-level Monte Carlo methods}


There is a close link between multi-fidelity methods and multi-level Monte Carlo methods.
Let us consider as control variates a \blue{hierarchy of discretizations} of the kinetic equation. For example, in the homogeneous case, with a cartesian grid we take 
\[
\Delta {\w}_h = 2^{1-h}(\Delta {\w}_1),\quad h=1,\ldots, L
\]
where $\Delta {\w}_1$ is the mesh width for the coarsest resolution, which corresponds to the solution with the lowest level of fidelity. Our full model is, therefore, represented by the fine scale solution obtained for $\Delta v_L$.
The hierarchy of numerical solutions $f_h(\theta,v,t)$, $h=1,\ldots, L$, at time $t$ with mesh $\Delta \w_h$ represents the setting for the multi-level control variate estimators. 

In particular, fixing all $\lambda_h=1$, $h=1,\ldots,L$, we get the classical \emph{Multi-level MC estimator} \cite{Giles}
\be
E_L^{{\bf 1}}[f](v,t) = E_{M_0}[f_1] + \sum_{h=1}^{L} (E_{M_h}[f_{h+1}-f_{h}]),
\label{eq:MLMC}
\ee
where we used the notation ${\bf 1}=(1,\ldots,1)^T$.  

Using the quasi-optimal values (or the optimal values) for $\lambda_h$ with the hierarchical grid constructed above we obtain a \green{quasi-optimal (optimal) MLMC} \cite{DPms2, HPW}. The main difference, compared to multi-fidelity models is the possibility to compute accuracy estimates between the various levels. On the other hand, the approach depends on additional parameters (the various grid sizes) which make its practical realization more involved to achieve optimal performances.

\begin{figure}[tb]
\begin{center}
\includegraphics[scale=0.25]{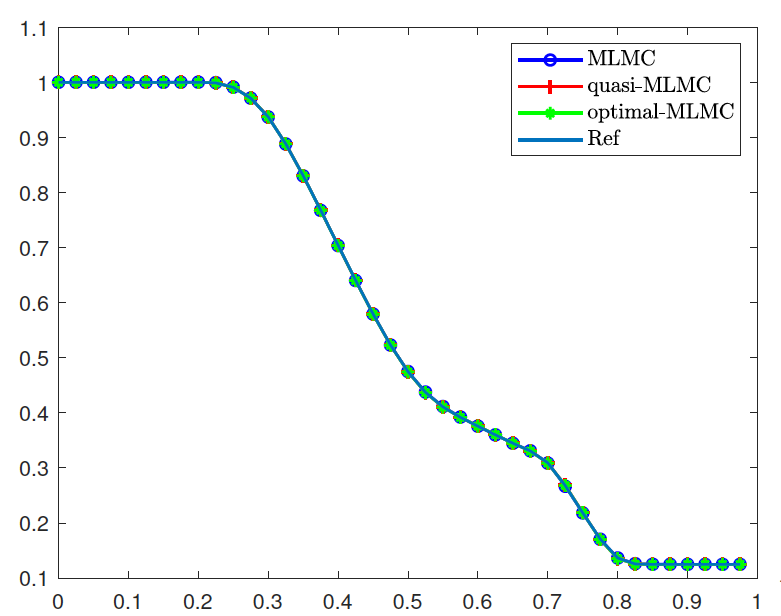}
\includegraphics[scale=0.25]{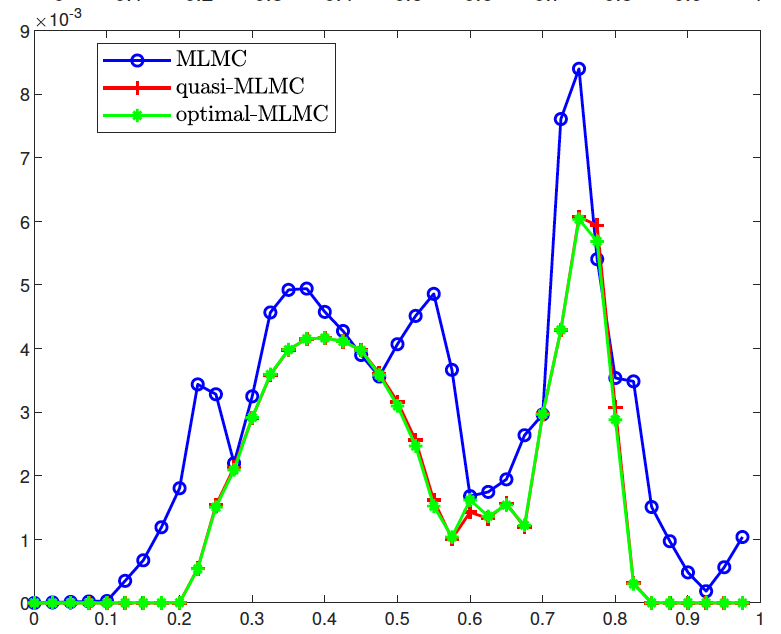}
\end{center}
\caption{Multilevel MC for BGK. $\E[\rho]$ (left) and error in space (right) for the various MLMC using $L=3$, $M_0 = 320$, $M_1 = 80$, $M_2 = 20$ with $\Delta x_h=2^{1-h}(\Delta x_1)$, $h=1,2,3$, $\Delta x_1=0.1$.}
\label{fg:bgkmlmc}
\end{figure}

In Figure \ref{fg:bgkmlmc} we report the results in the case of a space non homogeneous 
BGK model close to the fluid limit $\e=10^{-6}$ for the Sod shock tube problem \eqref{eq:sod} where $s=0.1$ and $\theta$ uniform in $[0,1]$. Improvements in the error curves obtained using the quasi-optimal and optimal MLMC over standard MLMC can be observed. 

\section{Structure preserving Stochastic-Galerkin (SG) methods}


For notation simplicity, let us assume a one-dimensional random space $\red{ z\in\Omega\subset \R}$, with $z$ distributed as $p(z)$, for a {space homogeneous} kinetic equation. We approximate $f(\theta,v,t)$ by its \blue{generalized Polynomial Chaos (gPC) expansion} \cite{Xu}
\be
f_M(\theta,v,t)=\sum_{m=0}^M \hat{f}_m(v,t)\Phi_m(\theta),
\label{eq:gPC}
\ee
where  $\left\{\Phi_{m}(\theta)\right\}_{m=0}^M$ are a set of orthogonal polynomials, of degree less or equal to $M$, {orthonormal} with respect to $p(z)$  
\[
\int_{\Omega} \Phi_n(\theta) \Phi_m(\theta) p(\theta)\,d\theta = \E[\Phi_m(\cdot)\Phi_n(\cdot)]=\delta_{mn},\qquad m,n=0,\ldots,M.
\] 
In \eqref{eq:gPC} the coefficients $\hat{f}_m$ are the {projection} of the solution with respect to $\Phi_m$
\be
\hat{f}_m(v,t)=\int_{\Omega} f(\theta,v,t)\Phi_m(\theta) p(\theta)\,d\theta=\E[f(\cdot,v,t)\Phi_m(\cdot)].
\ee
Stochastic Galerkin (SG) methods for kinetic equations based on the use of {deterministic methods} in the phase space have demonstrated numerical and theoretical evidence of {spectral accuracy}\cite{DJL, HJ, HJS, HJS2, JinPareschi, LJ, RHS}. However, their practical application presents some \blue{drawbacks}.


\begin{itemize}
\item 
SG methods lead to the \blue{loss of the physical properties}, like positivity and conservation of moments.
For example, if we denote by $m_\phi(f)=\int f \phi(v)\,dv$ the moments of $f$, we have 
\[
m_\phi(f) \neq m_\phi(f_M)= \sum_{m=0}^M m_\phi(\hat f_m)\Phi_m.
\]
These properties are essential to characterize the long time behavior of the system and the Maxwellian equilibrium states \cite{DP15, PZ1}.

\item One possibility, is to \blue{modify the coefficients} $\hat f_m$ in the gPC expansion or the \blue{polynomial basis} in such a way that the macroscopic moments of $f$ (or positivity) are preserved. This approach, however, is rather difficult in general and typically leads to the \blue{loss of spectral accuracy} \cite{CFY, Par3}.  

\item Additionally, for nonlinear hyperbolic conservation laws, like the Euler system in the fluid-dynamic limit, the generalized polynomial chaos expansion may lead to the \blue{loss of hyperbolicity} of the resulting approximated system (see \cite{DPe,DPL, HJS,JSh,PDL}). 
\end{itemize}


\subsection{Equilibrium preserving SG methods for the Boltzmann equation}
In this section we describe a general approach based on SG methods that permits to recover the correct long time behavior. To this aim,
let us consider the space homogeneous Boltzmann equation
\[
\partial_t f(\theta,v,t) = Q(f,f)(\theta,v,t).
\]
The \blue{standard SG method} reads
\[
\partial_t \hat f_h = \hat Q_h(\hat f,\hat f), \quad h=1,\ldots,M
\]
where $\hat f =(\hat f_0,\ldots,\hat f_M)^T$ and
\[
\hat Q_h(\hat f,\hat f) = \sum_{m,n=1}^M \hat f_m \hat f_n \E\left[Q(\Phi_m,\Phi_n)\Phi_h\right].
\]
Using the \green{decomposition} 
\[
f(\theta,{\w},t) = f^{\infty}(\theta,{\w})+g(\theta,{\w},t)
\]
from the bilinearity of $Q$ and the fact that $Q(f^\infty,f^\infty)=0$ we get
\[
Q(f,f)(\theta,{\w},t) = Q(g,g)(\theta,{\w},t) + \mathcal L(f^{\infty},g)(\theta,{\w},t),
\]
where $\mathcal L(\cdot,\cdot)$ is a linear operator defined as 
\[
\mathcal L(f^{\infty},g)(\theta,{\w},t) = Q(g,f^\infty)(\theta,{\w},t)+Q(f^\infty,g)(\theta,{\w},t).
\]



Thus we can apply the SG projection to the transformed problem
\[
\partial_t g(\theta,{\w},t) = Q(g,g)(\theta,{\w},t)+ \mathcal L(f^{\infty},g)(\theta,{\w},t),
\]
which admits $g^{\infty}(\theta,{\w})\equiv 0$ as unique equilibrium state. 

We can write the \emph{equilibrium preserving SG method} as
\begin{eqnarray*}
\partial_t \hat g_h &=& \hat Q_h(\hat g,\hat g) + \hat{\mathcal L}_h(\hat f^{\infty},\hat g),\quad h=0,\ldots,M\\
\hat f_h &=& \hat f_h^{\infty} + \hat g_h
\end{eqnarray*}
where $\hat g =(\hat g_0,\ldots,\hat g_M)^T$, $\hat f^{\infty} =(\hat f^\infty_0,\ldots,\hat f^\infty_M)^T$ and
\[
\hat{\mathcal L}_h(\hat f^{\infty},\hat g) = \sum_{m,n=1}^M \hat f^\infty_m \hat g_n \E\left[{\mathcal L}(\Phi_m,\Phi_n)\Phi_h\right].
\]
The values $\hat g_h=0$ (or equivalently $g^\infty_M=0$) are a local equilibrium of the SG scheme and
thus $\hat f_h = \hat f_h^{\infty}$ (or equivalently $f_M = f_M^{\infty}$) are a local equilibrium state. 


By substituting to $\hat g_h = \hat f_h - \hat f_h^{\infty}$ the SG scheme can be rewritten as
\[
\partial_t \hat f_h = \hat Q(\hat f,\hat f)-\hat Q(\hat f^\infty,\hat f^\infty), \quad h=1,\ldots,M,
\]
or equivalently 
\begin{eqnarray*}
\partial_t f_M(\theta,v,t) &=& Q_M(f_M,f_M)(\theta,v,t)-Q_M(f_M^\infty,f_M^\infty),\\
Q_M(f_M,f_M)(\theta,v,t)&=&\sum_{m=0}^M \hat Q_m(\hat f,\hat f)\Phi_m(\theta).
\end{eqnarray*} 
If we have a \blue{spectral estimate}%
for $Q(f,f)$, namely for $f\in H^r(\Omega)$ 
\[
\|Q(f,f)-Q_M(f_M,f_M)\|_{L^2(\Omega)} \leq \frac{C}{M^r}\left(\|f\|_{H^r(\Omega)}+\|Q(f_M,f_M)\|_{H^r(\Omega)}\right)
\]
and the equilibrium state $f^\infty\in H^r(\Omega)$, since $Q(f^\infty,f^\infty)=0$, we have
\[
\|Q_M(f^\infty_M,f^\infty_M)\|_{L^2(\Omega)} \leq \frac{C}{M^r}\left(\|f^\infty\|_{H^r(\Omega)}+\|Q(f^\infty_M,f^\infty_M)\|_{H^r(\Omega)}\right)
\]
which provide a spectral estimate for the equilibrium preserving SG method.

\subsection{Generalizations for nonlinear Fokker-Planck problems}

The approach just described applies to a large variety of kinetic equations where the equilibrium state is know. In the case of Fokker-Planck equations, the method can be generalized to the situation where the steady state is not known in advance \cite{DPZ2}. The idea is based on the notion of quasi-equilibrium state. To this aim given a one-dimensional Fokker-Planck equation characterized by
\begin{equation}
Q(f,f)=\partial_v \left(\mathcal{P}[f]f(\theta,v,t)+\partial_v (D(\theta,v)f(\theta,v,t))\right),
\end{equation}
we can consider solutions of the following problem
\begin{eqnarray*}
\mathcal{P}[f]f(\theta,v,t)+\partial_v (D(\theta,v)f(\theta,v,t))= 0,
\end{eqnarray*}
which gives
\begin{eqnarray*}
D(\theta,v)\partial_v f(\theta,v,t)=\left(\mathcal{P}[f](\theta,,v,t)+D'(\theta,v)\right)f(\theta,v,t).
\end{eqnarray*}
The above problem can be solved analytically for $f=f^\infty(v)$ only in some special cases. More in general we can represent a quasi-stationary solution in the form
\begin{equation}
f^q(\theta,v,t)=C\exp\left\{-\int_{-\infty}^v \frac{\mathcal{P}[f](\theta,v_*,t)+D'(\theta,v_*)}{D(\theta,v_*)}\,dv_*   \right\}
\end{equation}
being $C > 0$ a normalization constant. Therefore, $f^q$ is not the global in time equilibrium of the problem but have the property to annihilate the flux for each time $t \geq 0$ and that $f^q(\theta,v,t) \to f^{\infty}(v)$ as $t\to\infty$. 

Using the \green{decomposition} 
\[
f(\theta,{\w},t) = f^{q}(\theta,{\w})+g(\theta,{\w},t)
\]
it is clear that the formulation presented in Section 6.1 applies and we obtain a steady state preserving method for large times. We refer to \cite{DPZ2} for more details. 

\begin{figure}[tb]
\begin{center}
\includegraphics[scale=0.35]{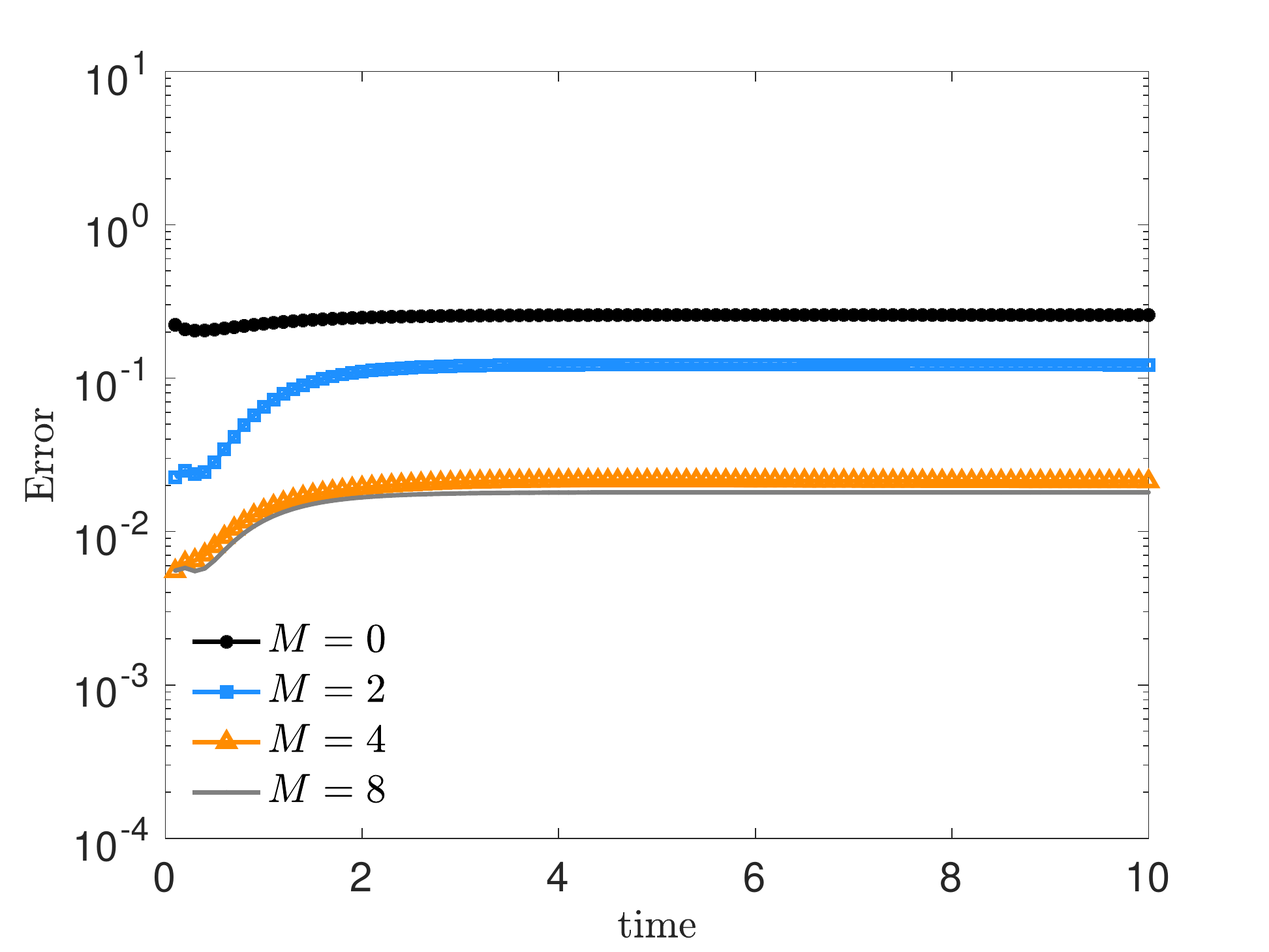}\hskip -.2cm
\includegraphics[scale=0.35]{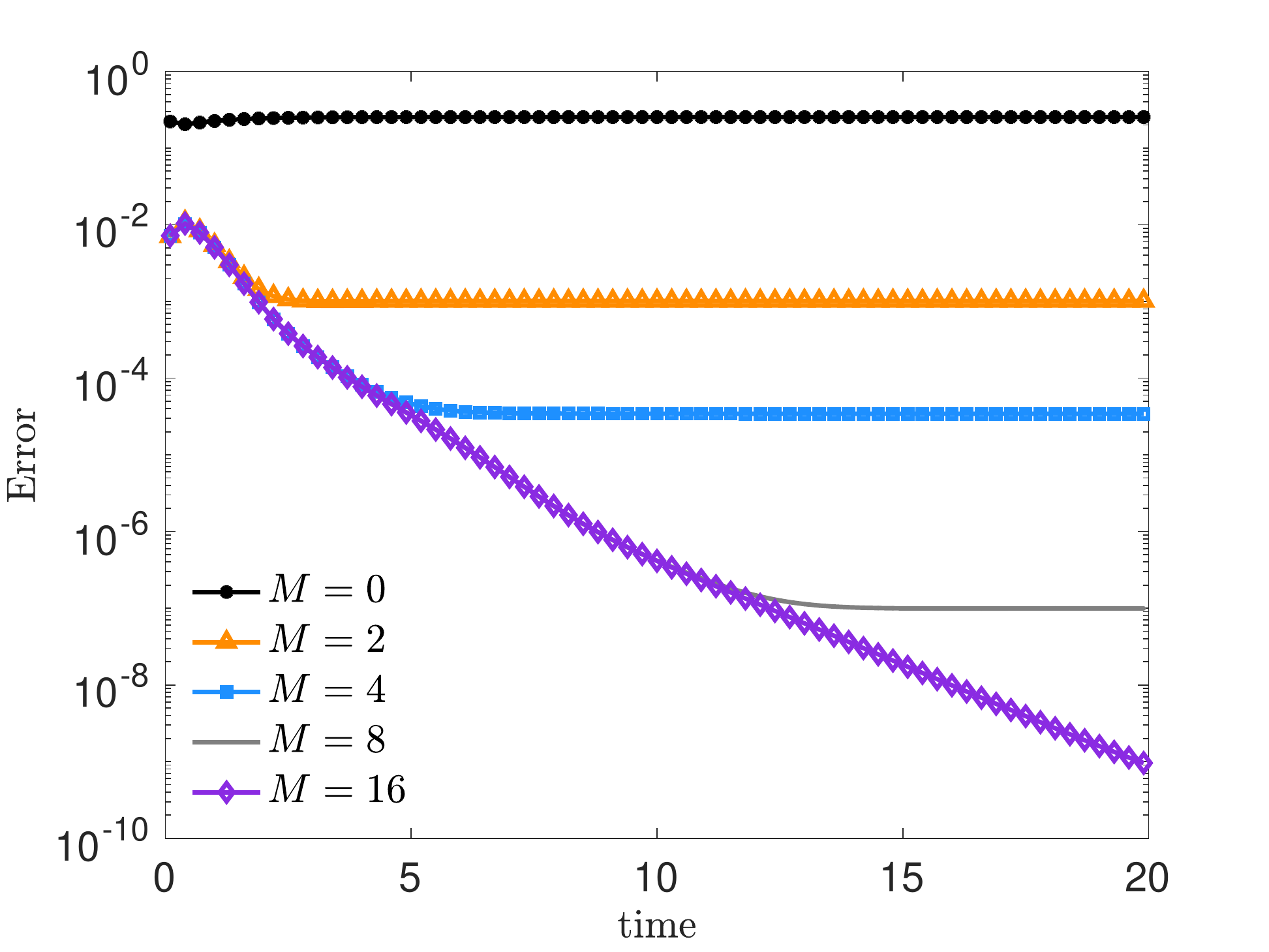}
\end{center}
\caption{Evolution of the $L^2$ error for the swarming model \eqref{eq:swarm} with standard SG scheme (left) and with the micro-macro SG scheme (right). The error has been computed with respect to a reference solution obtained with $M=40$, $N = 321$, $\Delta t = 10^{-1}$ and final time $T = 20$. }
\label{fig:swarm}
\end{figure}

As an example, let us consider the swarming model with self-propulsion defined by 
\begin{equation}
\mathcal{P}[f](\theta,v,t) = \alpha(\theta)(|v|^2-1)v + (v-u(\theta,t)),
\label{eq:swarm}
\end{equation}
where the mean velocity $u(\theta,t)=\int_{\RR^{d_v}} v f(\theta,v,t)\,dv$
is not conserved in time. In the above definition we assumed for simplicity $\int_{\RR^{d_v}}f(\theta,v,t)\,dv=1$. The quasi-stationary state is computed as
\begin{equation}
f^q(\theta,v,t)=C\exp\left\{-\frac1{D(\theta)}\left[\alpha(\theta)\frac{|v|^4}{4}+(1-\alpha(\theta))\frac{|v|^2}{2}-u_f(\theta,t)v\right]\right\}.
\end{equation}
In Figure \ref{fig:swarm} we report a comparison of the results obtained using a standard SG scheme and the micro-macro SG approach. We considered a diffusion coefficient $D(\theta) = 1/5 + \theta/10$ with $\theta\sim \mathcal U([-1,1])$ and deterministic self-propulsion $\alpha = 2$. In both cases the velocity space has been discretized by simple central differences using $N$ points in the domain $v\in[-2,2]$. It is evident that the error in the standard SG method saturates at the $O(\Delta v^2)$ order of the solver in the velocity space, while the micro-macro SG approach, thanks to its equilibrium preserving property, is able to achieve spectral accuracy asymptotically.


\section{Hybrid particle Monte Carlo SG methods}
The idea is to combine SG methods in the random space with particle Monte Carlo methods for the approximation of $f$ in the phase space. This  novel hybrid formulation makes it possible to construct efficient methods that preserve the main physical properties of the solution along with spectral accuracy in the random space \cite{CPZ, CZ, PZ2}. 

\subsection{Particle SG methods for Fokker-Planck equations}
We concentrate on a \blue{Vlasov Fokker-Planck} (VFP) for the evolution of $f = f(\theta,x,v,t)$ characterized by 
\[
Q(f,f)=\nabla_v \cdot \left(\mathcal P[f]f+\nabla_v (D f)\right)
\]
where
\[
\mathcal P[f](\theta,x,v,t) =  \int_{\mathbb R^{d_v}\times \mathbb{R}^{d_x}} P(\theta,x,x_*)(v-v_*)f(\theta,x_*,v_*,t)dv_* dx_*
\]
and diffusion $D(\theta)$.



The VFP equation can be derived from the following {system of stochastic differential equations} for $(X_i(\theta,t),V_i(\theta,t)) \in \mathbb R^{d_v}\times \mathbb{R}^{d_x}$, $i = 1,\dots,N$ with random inputs 
\[
\begin{cases}
dX_i(\theta,t) = V_i(\theta,t)dt \\
dV_i(\theta,t) = \dfrac{1}{N} \displaystyle\sum_{j=1}^N P(\theta,X_i,X_j)(V_j-V_i)dt + \sqrt{2D(\theta)}dW_i,
\end{cases}
\]
being $\{W_i\}_{i=1}^N $ independent Brownian motions. 

We consider the \emph{empirical measure} associated to the particle system 
\[
f^{(N)}(\theta,x,v,t) = \dfrac{1}{N} \sum_{i = 1}^N \delta(x-X_i(\theta,t))\otimes \delta(v-V_i(\theta,t))
\]
Under suitable assumptions it can be shown that as $N\to\infty$, the empirical measure $f^{(N)}\to f$ solution of the VFP problem \cite{CPZ}.



We consider the SG approximation of the particle system, given by
\[
X_i^M (\theta,t)= \sum_{m = 0}^M \hat{X}_{i,m}(t)\Phi_m(\theta),\qquad V_i^M(\theta,t) = \sum_{m=0}^M \hat{V}_{i,m}(t)\Phi_m(\theta),\quad i=1,\ldots,N
\]
where $\hat{X}_{i,m}$, $\hat{V}_{i,m}$ are the {projections} of the solution with respect to $\Phi_m$
\[
\hat{X}_{i,m}(t)=\E[X_i(\cdot,t)\Phi_m(\cdot)],\qquad \hat{V}_{i,m}(t)=\E[V_i(\cdot,t)\Phi_m(\cdot)].
\]
The \blue{particle SG method} is then obtained as
\[
\begin{cases}
d\hat{X}_{i,h} &= \hat{V}_{i,h}dt \\[+.2cm]
d\hat{V}_{i,h} &= \displaystyle\dfrac{1}{N}\sum_{j = 1}^N \sum_{k = 0}^M P_{hk}^{ij} (\hat{V}_{j,k}-\hat{V}_{i,k})dt + D_h dW_i
\end{cases}
\]
and $P_{hk}^{ij} = \mathbb E[P(\cdot,X_i^M,X_j^M) \Phi_h(\cdot)\Phi_k(\cdot)]$, $D_h = \mathbb E[ \sqrt{2 D(\cdot)}\Phi_h(\cdot)]$.

Moments are recovered from the empirical measure as
\begin{eqnarray*}
f_M^{(N)}(\theta,x,v,t) &=& \dfrac{1}{N} \sum_{i = 1}^N \delta(x-X^M_i(\theta,t))\otimes \delta(v-V^M_i(\theta,t))\\
m_\phi(f_M^{(N)})&=&\dfrac{1}{N} \sum_{i = 1}^N \delta(x-X^M_i(\theta,t))\phi(V^M_i(\theta,t))
\end{eqnarray*}


%

The method just described has the usual \blue{quadratic cost} $O(N^2)$ of a mean field problem, where each particle at each time step modifies its velocity interacting with all other particles. In addition, this cost has to be multiplied by the quadratic cost $O(M^2)$ of the SG method. Therefore the overall computational cost is $O(M^2N^2)$.
 
 A reduction of the cost is obtained using a suitable \green{Monte Carlo evaluation} of the interaction dynamics 
 to mitigate \blue{the curse of dimensionality} \cite{albipar13} 
\[
\begin{cases}
d\hat{X}_{i,h} &= \hat{V}_{i,h}dt \\[+.2cm]
d\hat{V}_{i,h} &= \displaystyle\red{\dfrac{1}{S}\sum_{j \in S_i}} \sum_{k = 0}^M P_{hk}^{ij} (\hat{V}_{j,k}-\hat{V}_{i,k})dt + D_h dW_i
\end{cases}
\]
where $S_i$ is a random subset of size $S\leq N$ of the particles indexes $\{1,2,\ldots,N\}$.


Using the Euler-Maruyama method to update the particles we have the following algorithm.

\begin{algorithm}[Particle SG algorithm]~
{\sl 
\begin{enumerate}
\item Consider $N$ samples $(X_i,V_i)$ from $f_0(x,v) $ and fix $S \le N$.
\item Perform gPC representation on the particles: $\left(\hat X_{i,h},\hat V_{i,h} \right)$, for $h = 0,\dots,M$
\item For $n = 0,\dots,T-1$ 
\begin{itemize}
\item Generate $N$ random variables $\{\eta_i\}_{i=1}^N \sim \mathcal N(0,1)$
\item For $i = 1,\dots,N$
\begin{itemize}
\item Sample $S$ particles $\{j_1,\dots,j_{S}\} := S_i$ uniformly without repetition
\item Compute the space and velocity change
\[
\begin{split}
\hat X_{i,h}^{n+1} &= \hat X_{i,h}^n + \hat V_{i,h}^n \Delta t \\
\hat V_{i,h}^{n+1}  &= \hat V_{i,h}^n + \dfrac{\Delta t}{S} \sum_{j \in S_i} \sum_{k = 0}^M P_{kh}^{i j} (\hat V_{j,k}^n-\hat V_{i,k}^n) + \sqrt{\Delta t}\,D_h \, \eta_i
\end{split}
\]
\end{itemize}
\end{itemize}
\item Reconstruct the quantity of interest $\mathbb E[\Psi(F(f))]$
\end{enumerate}
}
\end{algorithm}
The last step can be performed directly using the empirical distribution $f^{(N)}$ or some suitable reconstruction of $f$ using standard techniques.
Thanks to the random subset evaluation of the interaction sum the overall \blue{cost is reduced to $O(M^2 S N)$}, with $S \ll N$.

\begin{figure}[tb]
\begin{center}
\includegraphics[scale=0.5]{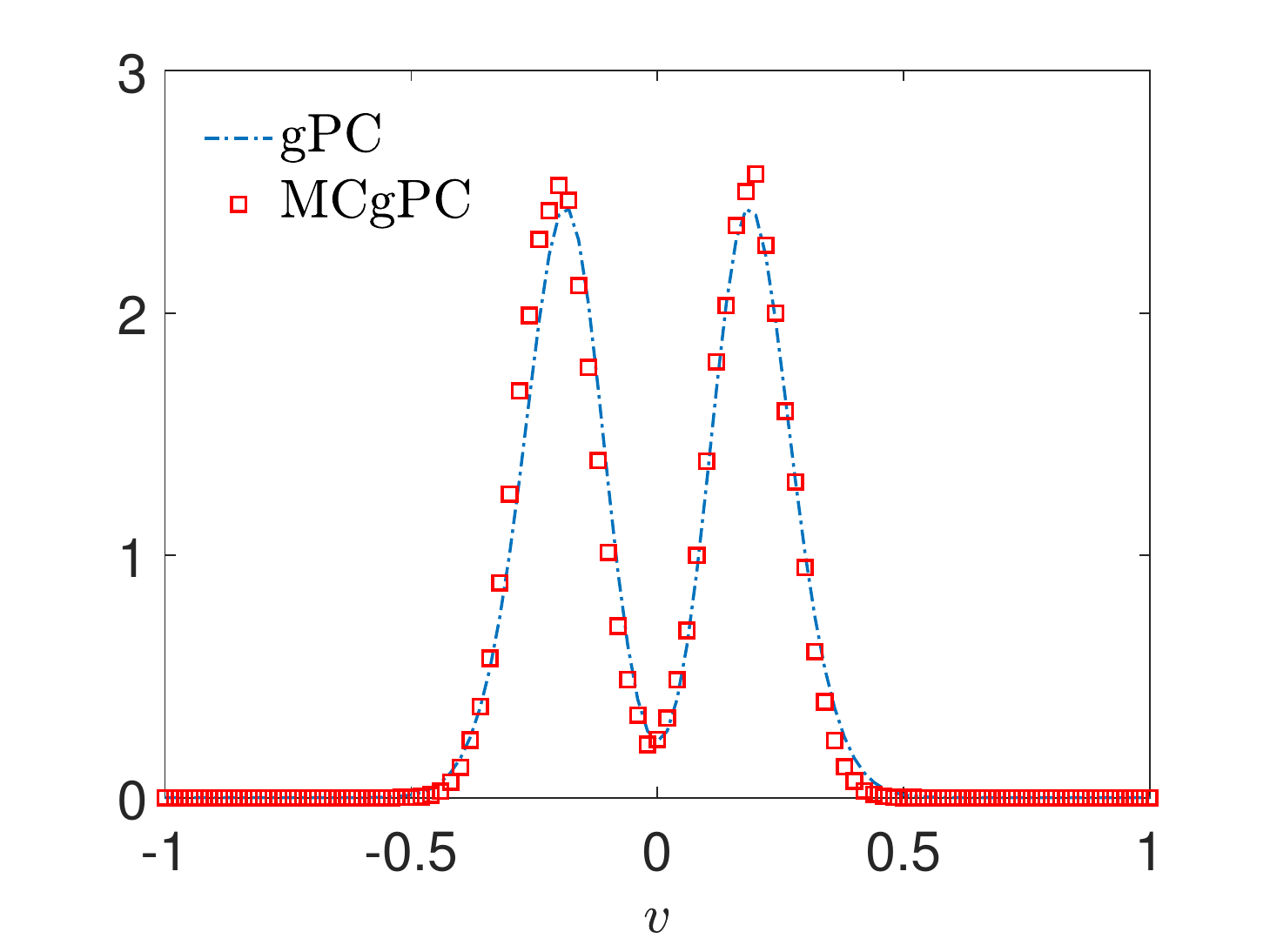}
\includegraphics[scale=0.5]{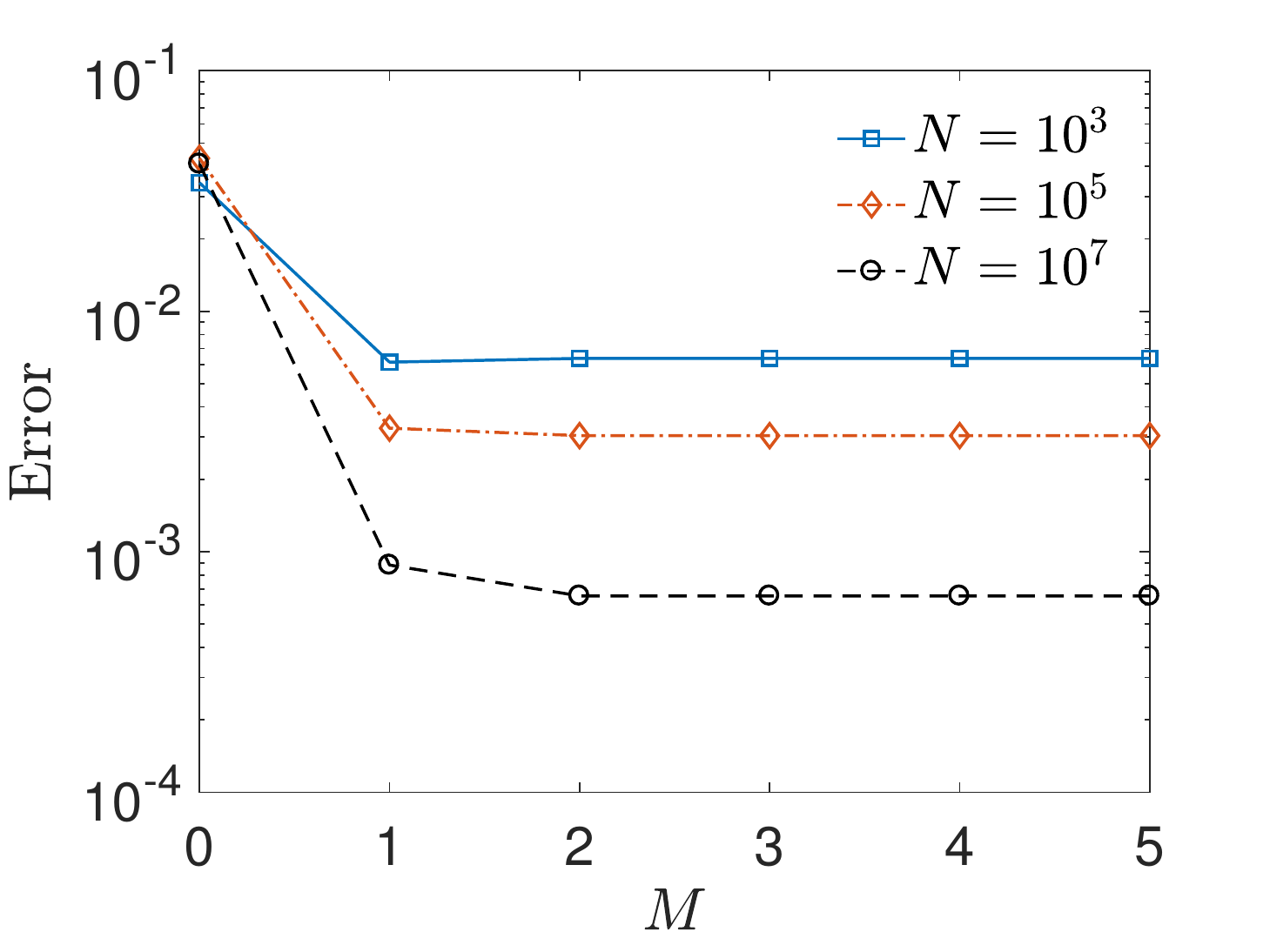}
\end{center}
\caption{Left: expected density at time $T = 1$ obtained through a standard SG method (gPC) and the particle gPC scheme (MCgPC) with $S = 5$ at each time step. The gPC expansion has been performed up to order $M = 5$. Right: convergence for the expected temperature of the MCgPC method with fixed $S = 50$ and an increasing number of particles. The reference temperature has been computed with a standard SG method for the mean-field problem. }
\label{fig:PgPC}
\end{figure}

\begin{remark}~
\begin{itemize}

\item The advantage of considering a SG scheme for the particle system lies in the preservation of the typical \blue{spectral convergence} in the random space together with the \blue{physical properties} of the original system.  

\item In the case $S = N$ we obtain the typical convergence rate $O(N^{-1/2})$ due to Monte Carlo sampling in the phase space. The fast evaluation of the interactions induces an additional error $O\left(\sqrt{{1}/{S}-{1}/{N}}\right)$ with $S<N$. 

\end{itemize}
\end{remark}

We report in Figure \ref{fig:PgPC} the result of a simulation concerning the simple space homogeneous one-dimensional alignment process corresponding to $P(\theta,x,x_*)=1+s\theta$, $\theta \sim U(0,1)$, $s=0.5$ and $D=0$. The initial data is given by a bimodal density
\[
f_0(v)=\beta\left[e^{\frac{(v-\mu)^2}{2\sigma^2}}+e^{\frac{(v+\mu)^2}{2\sigma^2}}\right]
\]
with $\sigma^2=0.1$, $\mu=0.25$ and $\beta$ such that $\int_{\RR} f_0(v)\,dv=1$. It is clear that a very small value $M$ suffices to match the accuracy in the random and in the phase space. 



\subsection{Direct simulation Monte Carlo SG methods}
The extension of the particle SG approach just discussed to Boltzmann type equations is non trivial. We recall here the basic methodology in the simple case of Maxwell molecules and refer to \cite{PZ2} for details on its extension to the variable hard sphere cases.

To this aim, we will focus on the space homogeneous Boltzmann equation, and observe that, in the case of Maxwell molecules $B\equiv 1$, the collision operator can be rewritten as
\begin{equation}
Q(f,f)(z,v,t)=Q^+(f,f)(z,v,t)-\mu f(z,v,t),
\label{eq:BMMi}
\end{equation}
where $\mu>0$ is a constant and we assumed $\int_{\R^{d_v}} f(z,v_*,t)\,dv_*=1$, $\forall\,z \in \Omega$.

We consider a set of $N$ samples $v_i(z,t)$, $i=1,\ldots,N$ from the kinetic solution at time $t$ and approximate $v_i(z,t)$ by its {generalized polynomial chaos (gPC) expansion}
\[
v^M_i(z,t)=\sum_{m=0}^M \hat{v}_{i,m}(t)\Phi_m(z).
\]
where  $\hat{v}_{i,m}(t)=\E[v_i(\cdot,t)\Phi_m(\cdot)]$.

To define the DSMC-SG algorithm we consider the projection on the above space of the collision process in the DSMC method (see \cite{Par2}). In the case of the uncertain Boltzmann collision term \eqref{eq:BMMi} we have
\begin{eqnarray*}
v_i'(z,t) &=& \frac12(v_i(z,t)+v_j(z,t))+\frac12 |v_i(z,t)-v_j(z,t)|\omega,\\
v_j'(z,t) &=& \frac12(v_i(z,t)+v_j(z,t))-\frac12 |v_i(z,t)-v_j(z,t)|\omega.
\end{eqnarray*}
Let us observe that
\be
|v_i'(z,t)-v_j'(z,t)| =  |v_i(z,t)-v_j(z,t)|,
\label{eq:mrv}
\ee
so that the modulus of the relative velocity is unchanged during collisions. 

We first substitute the velocities by their gPC expansion   
\begin{eqnarray*}
{v^{M}_i}^\prime(z,t) &=& \frac12(v^M_i(z,t)+v^M_j(z,t))+\frac12 |v^M_i(z,t)-v^M_j(z,t)|\,\omega,\\
{v^{M}_i}^\prime(z,t) &=& \frac12(v^M_i(z,t)+v^M_j(z,t))-\frac12 |v^M_i(z,t)-v^M_j(z,t)|\,\omega
\end{eqnarray*}
and then project by integrating against $\Phi_m(z)\,p(z)$ on $\Omega$ to get for $m=0,\ldots,M$
\begin{eqnarray}
\label{eq:c1}
\hat{v}_{i,m}^{\prime}(t) &=& \frac12(\hat{v}_{i,m}(t)+\hat{v}_{j,m}(t))+\frac12 \hat V^m_{ij}\,\omega,\\
\label{eq:c2}
\hat{v}_{j,m}^{\prime}(t) &=& \frac12(\hat{v}_{i,m}(t)+\hat{v}_{j,m}(t))-\frac12 \hat V^m_{ij}\,\omega
\end{eqnarray}
where
\be
\hat V^m_{ij}=\int_{\Omega} |v^M_i(z,t)-v^M_j(z,t)|\Phi_m(z)\, p(z)\,dz,
\label{eq:cmt}
\ee
is a time independent matrix thanks to \eqref{eq:mrv} which consists of a total of $(M+1)N^2$ elements that can be computed accurately and stored once for all at the beginning of the simulation.

\begin{figure}[tb]
\centering
\subfigure[$t = 0$]{
\includegraphics[scale = 0.25]{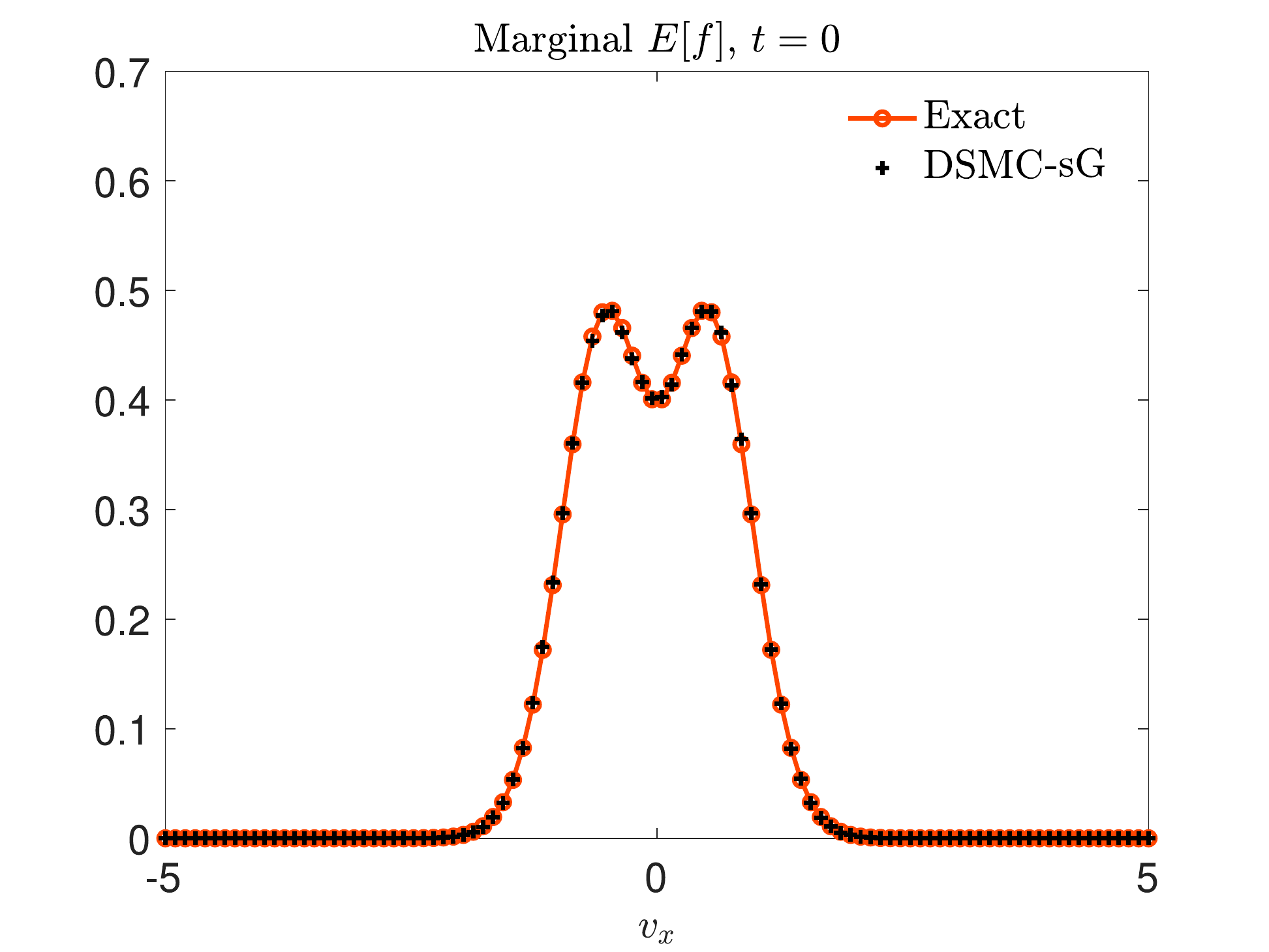}}\hskip -.2cm
\subfigure[$t = 1$]{
\includegraphics[scale = 0.25]{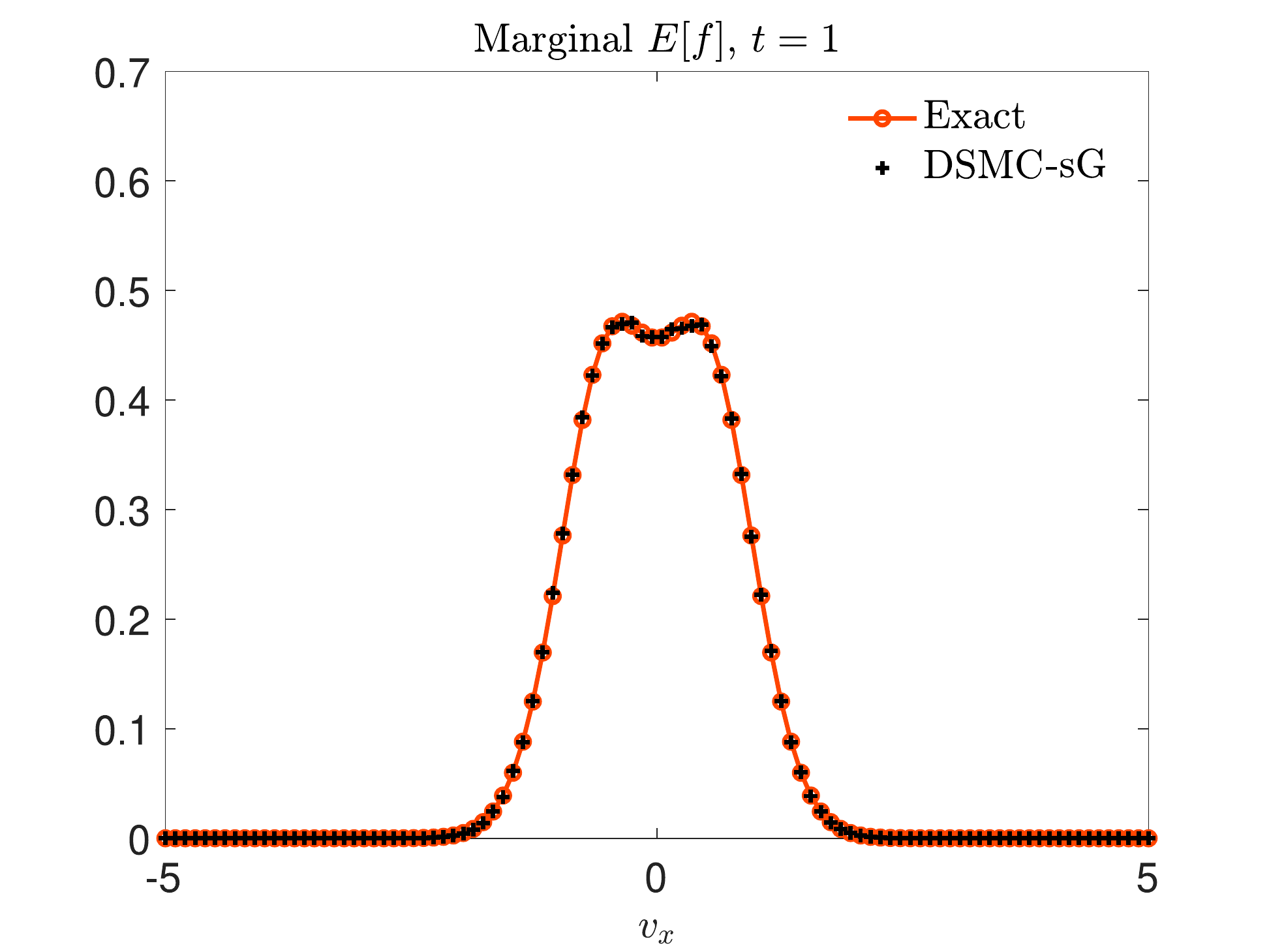}}\hskip -.2cm
\subfigure[$t = 5$]{
\includegraphics[scale = 0.25]{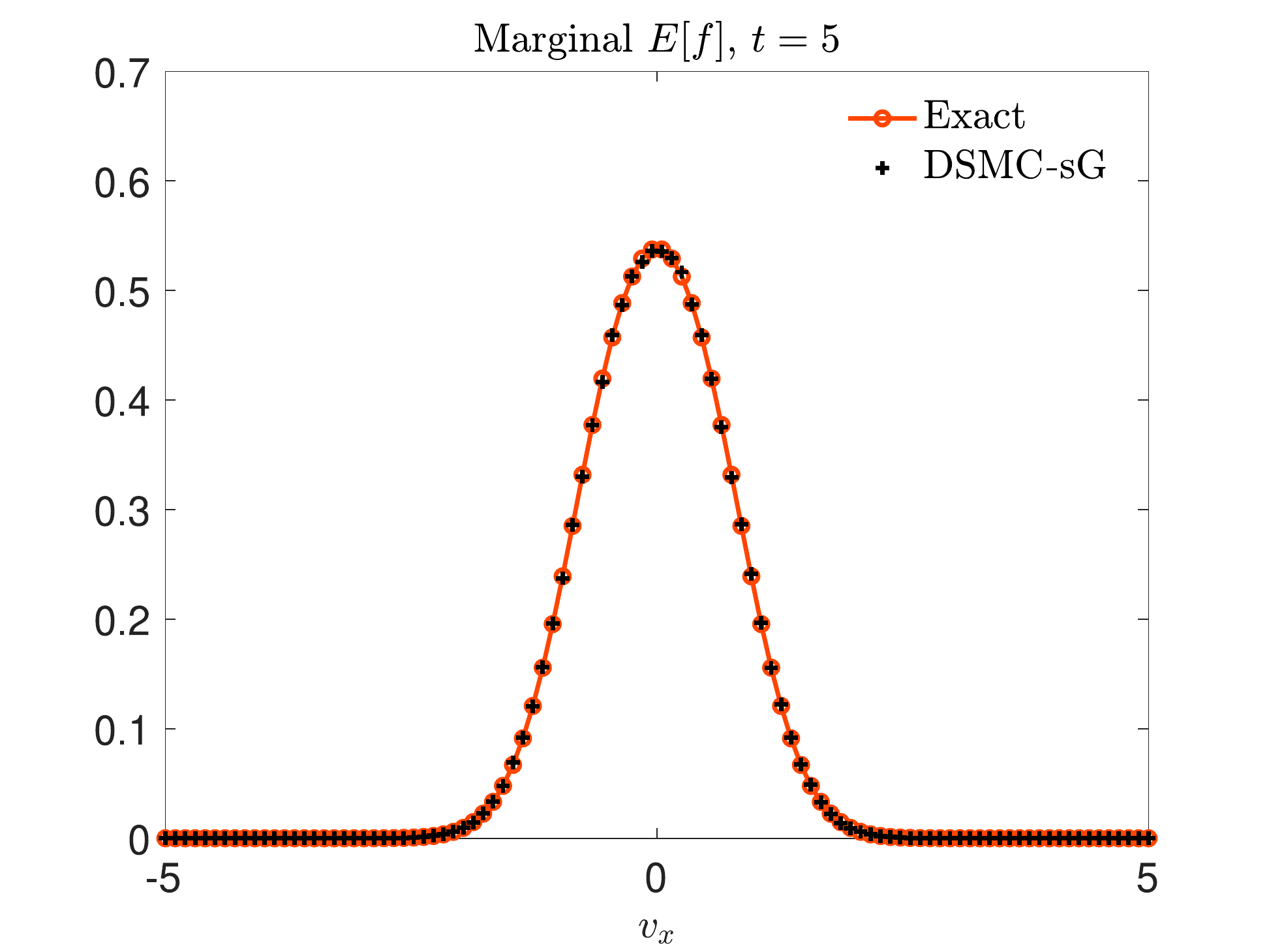}} \\
\subfigure[$t = 0$]{
\includegraphics[scale = 0.25]{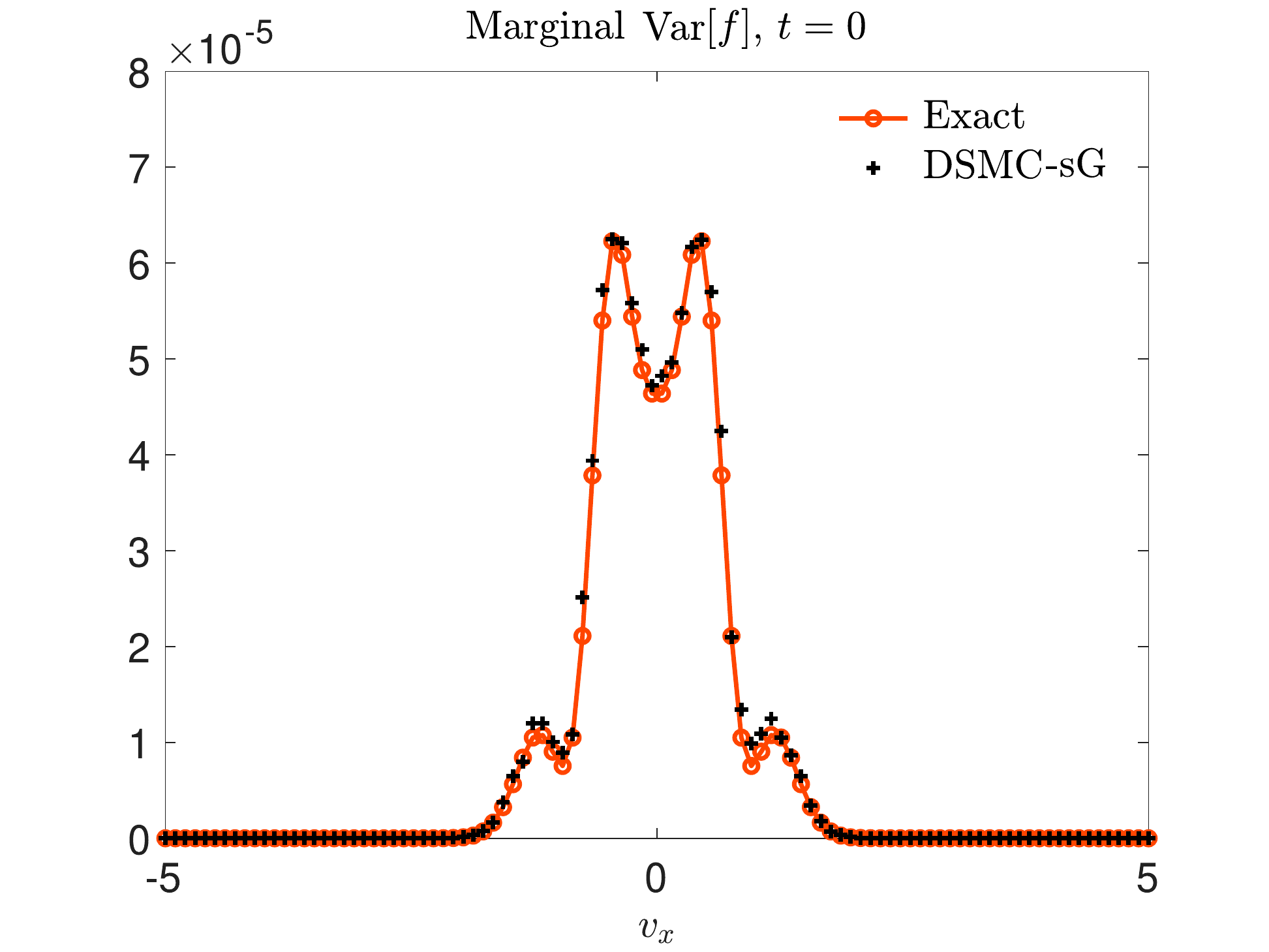}}\hskip -.2cm
\subfigure[$t = 1$]{
\includegraphics[scale = 0.25]{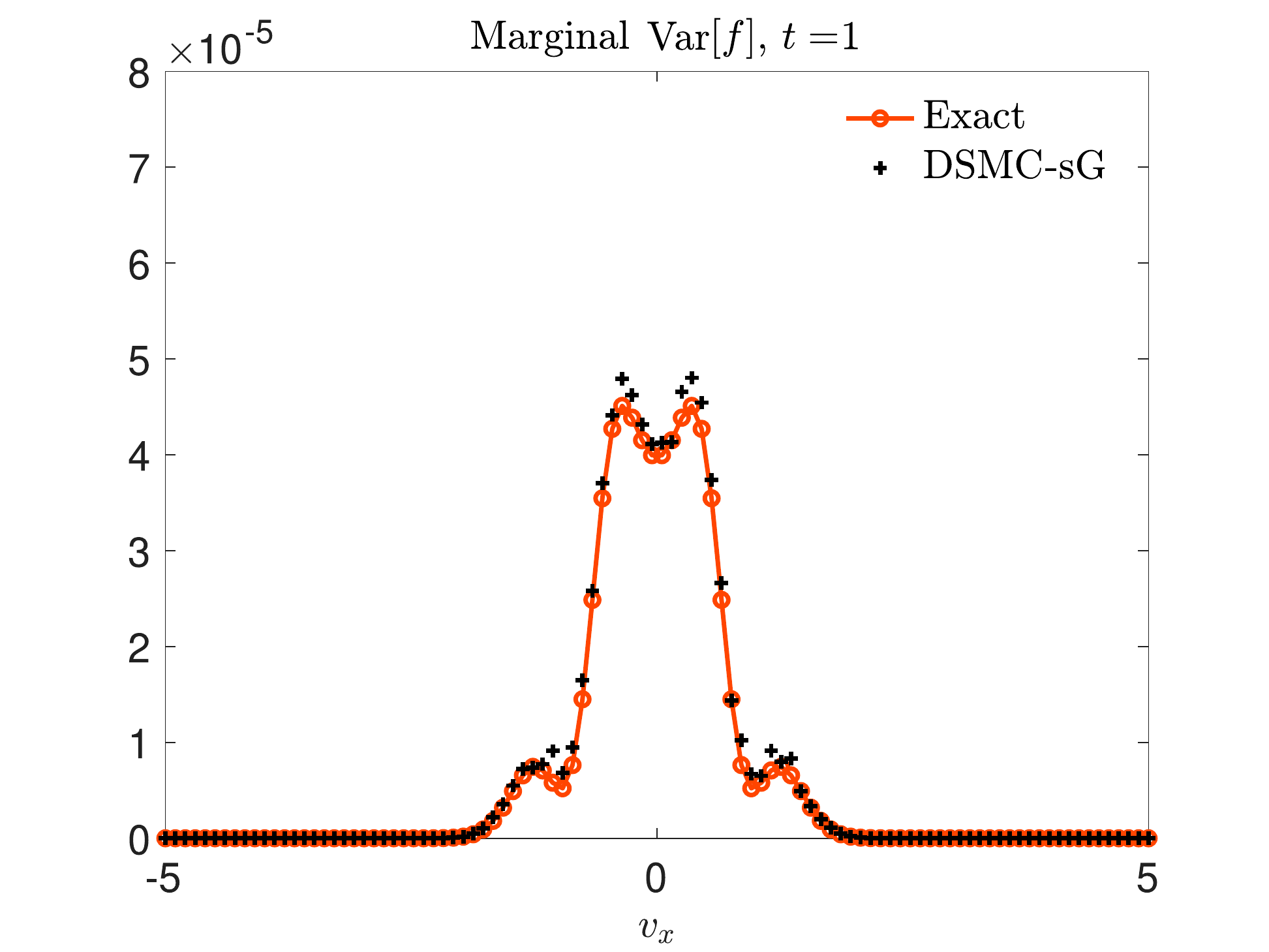}}\hskip -.2cm
\subfigure[$t = 5$]{
\includegraphics[scale = 0.25]{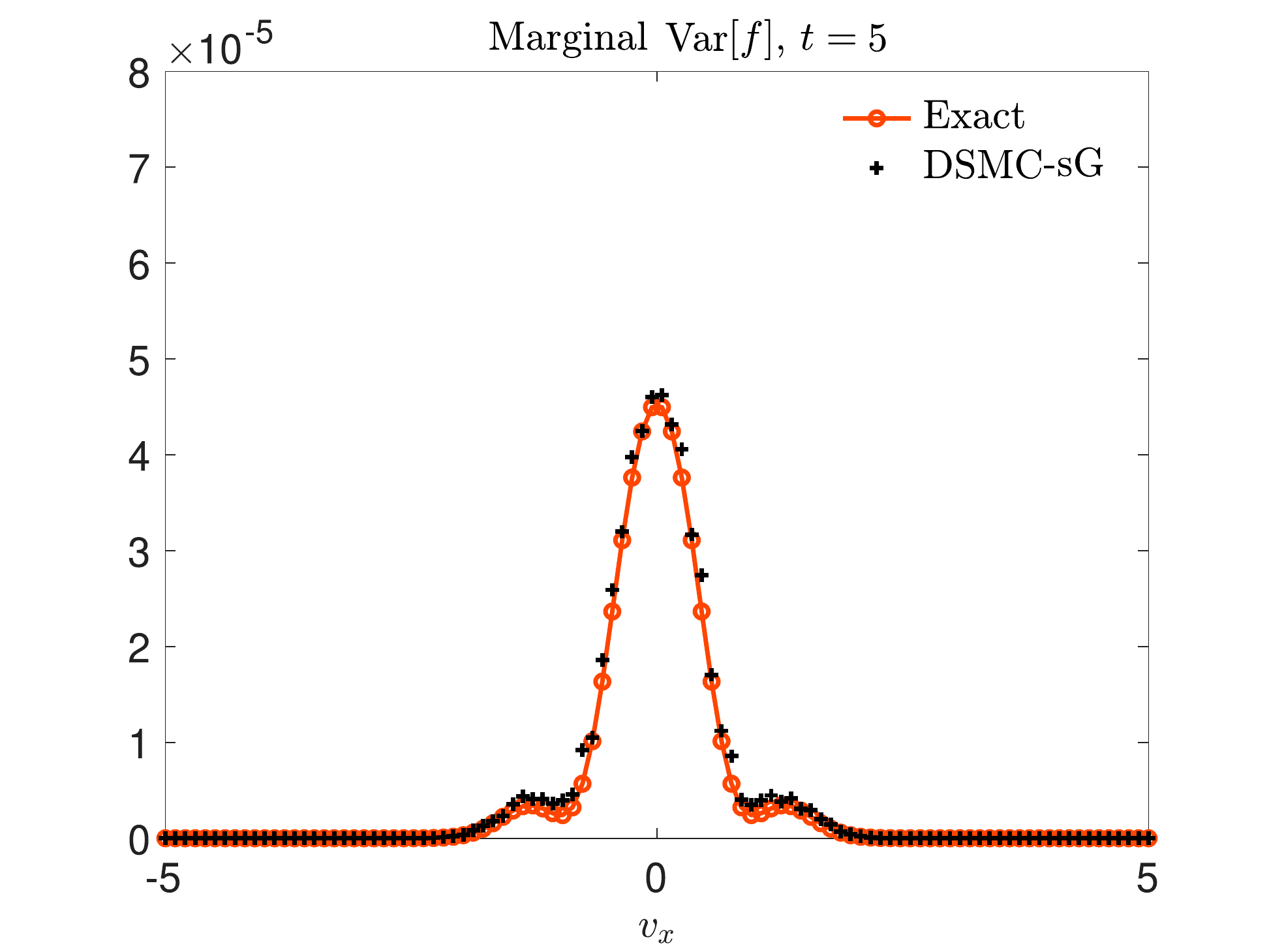}}
\caption{Evolution at times $t = 0,1,5$ of the marginal $\mathbb E[f]$ and $\textrm{Var}(f)$ from exact solution \eqref{eq:exact_max} and DSMC-SG approximation of the 2D Boltzmann model for Maxwell molecules with uncertain temperature. We considered $N = 10^6$ particles with $M= 5$ Galerkin projections and $\Delta t = 10^{-1}$. The reconstruction step has been performed in $[-5,5]^2$ through $100^2$ gridpoints. }
\label{fig:max3}
\end{figure}

Thus, each Monte Carlo collision can be performed at a computational cost of $O(M)$ which is the minimum cost to update the $M$ modes of each velocity. The SG extensions of the DSMC algorithms by Nanbu for Maxwell molecules is reported below.

\begin{algorithm}[DSMC-SG for Maxwell molecules]~
\begin{enumerate}
\item Compute the initial gPC expansions $ \{v^{M,0}_i, i=1,\ldots,N\} $,\\ 
      from the initial density $f_0(v)$
\item Compute the collision matrix $\hat V^m_{ij}$, $i,j=1,\ldots,N$, $m=0,\ldots,M$,\\
using \eqref{eq:cmt}. 
\item   \begin{tabbing}
\= {\fP for} \= $n=0$  {\fP to} $\nt-1$ \\
          \>          \> given $\{\hat v_i^{m,n},i=1,\ldots,N,\,m=0,\ldots,M\}$\\
          \>       \>   \ind \= \cir \= set $N_c = \IR(\mu N\Delta t/2)$ \\
          \>       \>        \> \cir \> select $N_c$ pairs $(i,j)$ uniformly among all possible pairs, \\
          \>       \>        \>      \> - \= perform the collision between $i$ and $j$, and compute \\
          \>       \>        \>      \>      \> $\hat v_{i,m}^{\prime}$ and $\hat v_{j,m}^{\prime}$ according to \eqref{eq:c1}-\eqref{eq:c2} \\
          \>       \>        \>      \> - set $\hat v_{i,m}^{n+1} = \hat v_{i,m}^{\prime}$, $\hat v_{j,m}^{n+1}=\hat v_{j,m}^{\prime}$ \\
          \>       \>        \> \cir \> set $\hat v_{i,m}^{n+1}=\hat v_{i,m}^{n}$ for all the particles that have not been selected\\
         \> {\fP end for}
\end{tabbing}
\end{enumerate}
\end{algorithm}

As a numerical example we consider the 2D case with uncertain initial data corresponding to the exact solution\cite{Boby, PZ2} 
\begin{equation}\label{eq:exact_max}
f(z,v,t) = \dfrac{1}{2\pi s(z,t)} \left[ 1 - \dfrac{1-\alpha(z)s(z,t)}{\alpha(z)s(z,t)} \left(1-\dfrac{\mathbf v^2}{2s(z,t)}  \right)\right]e^{-\frac{\mathbf v^2}{2s(z,t)}},
\end{equation}
where $s(z,t) = \dfrac{2 -e^{-t/8}}{2\alpha(z)}$. We will consider $\alpha(z) = 2+\kappa z$, with $z \sim  U(-1,1)$.

\begin{figure}[tb]
\centering
\includegraphics[scale = 0.35]{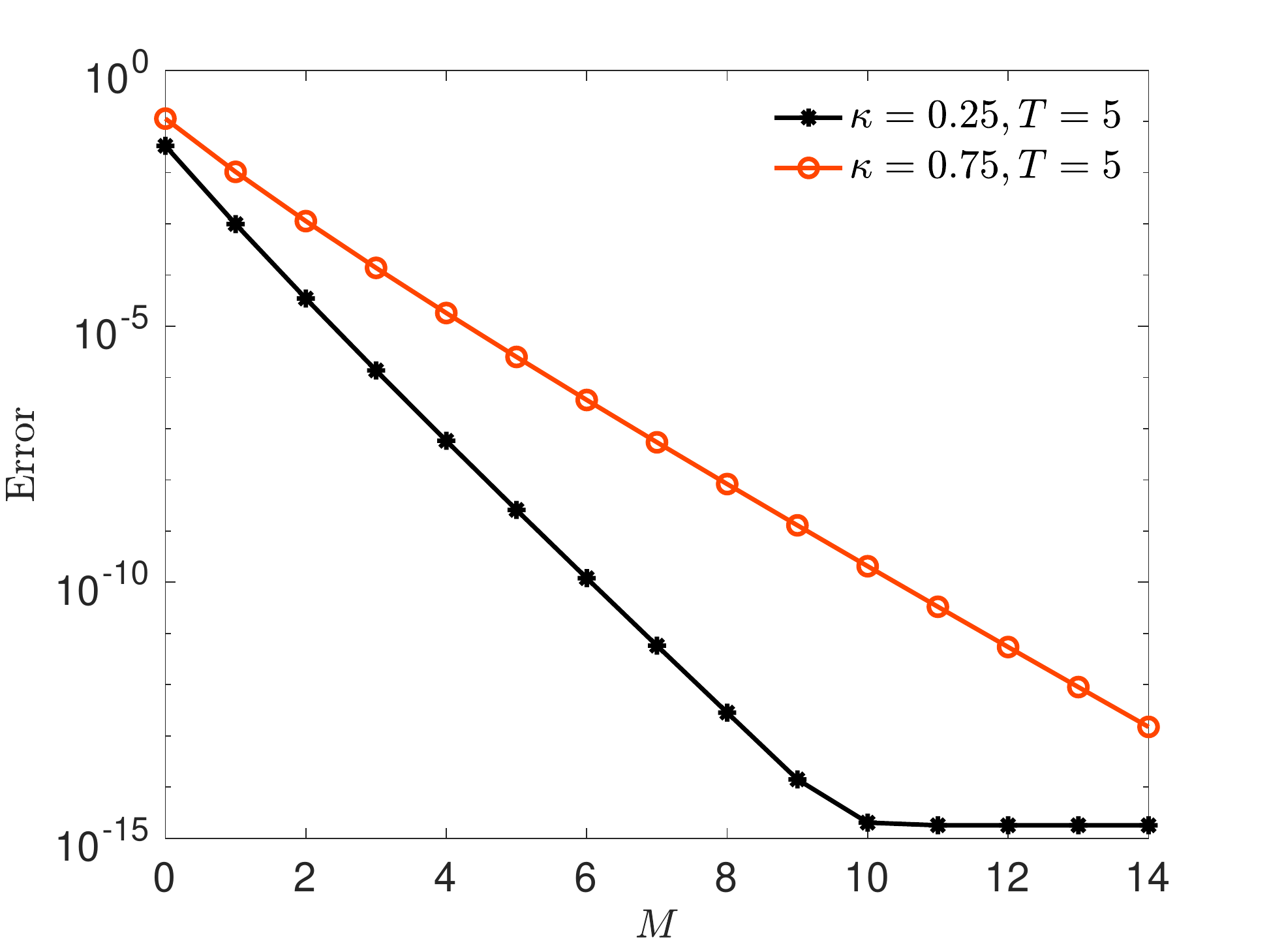}
\includegraphics[scale = 0.35]{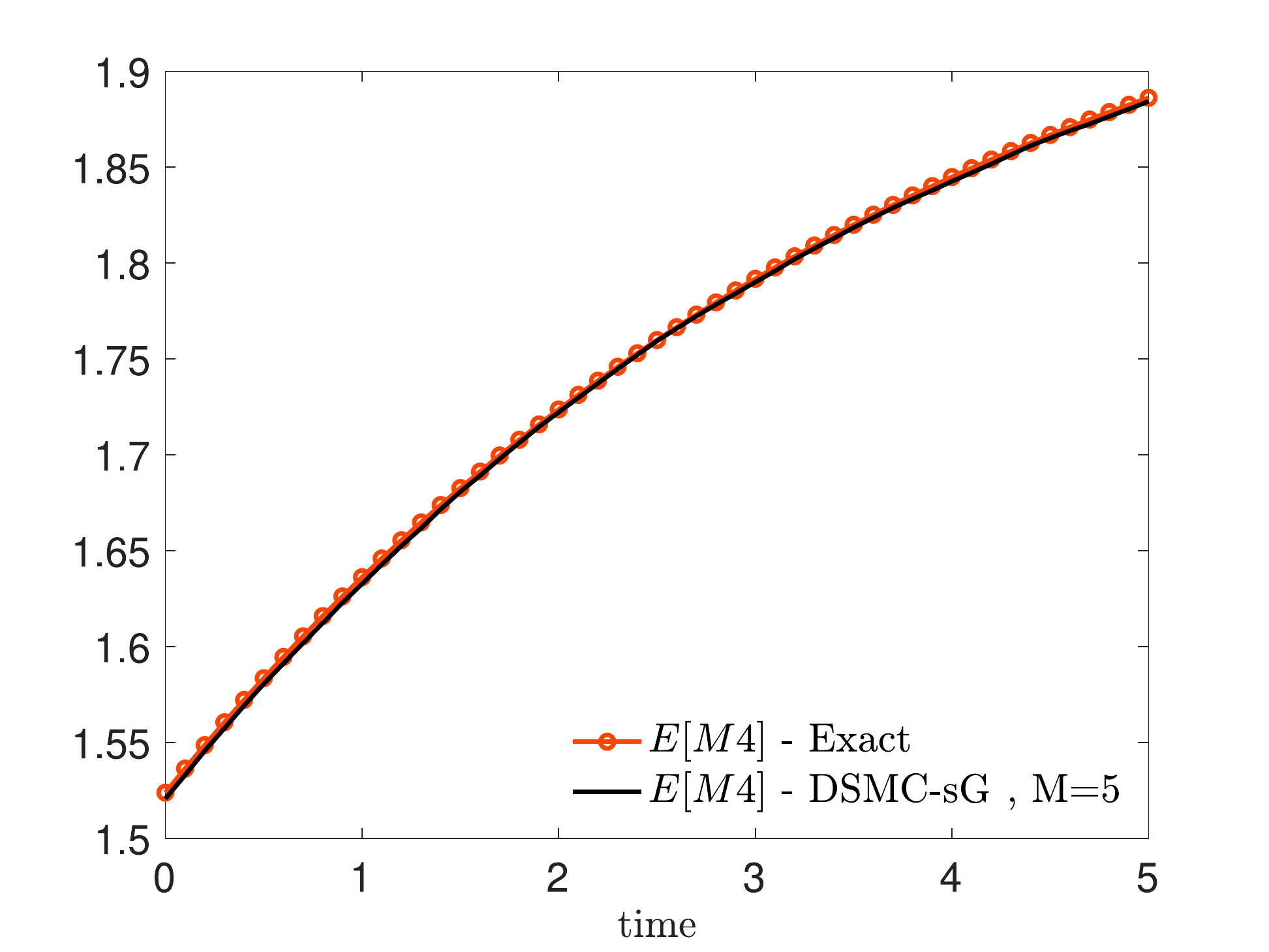}
\caption{Left: Convergence of the $L^2(\Omega)$ error with respect to the fourth order moment obtained from a reference solution computed with $N = 10^6$ particles and $M = 25$ from the DSMC-SG methods. Right: evolution of the fourth order moment in the interval $[0,5]$ for exact and DSMC-SG approximation with $N = 10^6$ and $M = 5$.  }
\label{fig:max4}
\end{figure}

To emphasize the good agreement of the computed approximation for all times, we depict in Figure \ref{fig:max3} the evolution at times $t = 0,1,5$ of the marginal of $\mathbb E[f]$ and $\textrm{Var}(f)$. 
Next, in Figure \ref{fig:max4} we present spectral convergence of the scheme computed through the fourth order moment of the 2D model with $\alpha(z) = 2+\kappa z$, $\kappa = 0.25$ and $\kappa = 0.7 5$ with $z\sim  U(-1,1)$. As reference solution we considered the fourth order moment at time $T = 5$ obtained with $N = 10^6$ particles and $M = 25$ Galerkin projections and the evolution is computed with $\Delta t = 10^{-1}$. In the right plot we present the decay of the $L^2(\Omega)$ error for increasing $M = 0,\dots,14$ in semilogarithmic scale. In the left plot we represent also the whole evolution of $M4$ computed through exact solution and through its DSMC-SG approximation. We obtain numerical evidence of spectral convergence.

\begin{remark}~
In a space non-homogeneous setting, the relative velocity changes at each time step due to the transport process.  Therefore, the largest part of the computational cost is due to the computation of the matrix \eqref{eq:cmt} at each time step. Note, however, that since $i$ and $j$ are selected at random, we may not need all elements in the matrix in the collision process. Thus, for fixed values of $i$ and $j$ we approximate the vector $\hat V^m_{ij}$ by Gauss quadrature 
\be
\hat V^m_{ij}(t) \approx \sum_{h=0}^H w_h |v^M_i(z_h,t)-v^M_j(z_h,t)|\Phi_m(z_h).
\ee 
The resulting scheme requires $O(MH)$ operations to compute $v^M_i(z_h,t)$ and $v^M_j(z_h,t)$ for all $h$'s and $O(MH)$ operations to evaluate $\hat V^m_{ij}(t)$ for all $m$'s. Taking $H=M$ the total cost of a Monte Carlo collision at each time step is therefore $O(M^2)$. 
\end{remark}

\subsection*{Acknowledgments} This work has been supported by the Italian Ministry of Instruction, University and Research (MIUR) under the PRIN Project 2017, No. 2017KKJP4X, "Innovative numerical methods for evolutionary partial differential equations and applications".

%
%
%

\end{document}